\documentclass[12pt]{amsart}

\usepackage{amsmath,amsfonts,amssymb,graphicx,hyperref}
\usepackage{amsthm}

\newcommand{\paren}[1]{\ensuremath{\left( #1 \right)}}
\newcommand{\set}[1]{\ensuremath{\left\{ #1 \right\}}}
\newcommand{\abs}[1]{\ensuremath{\left| #1 \right|}}
\newcommand{\setdiv}{\,\middle|\,}
\newcommand{\summod}[1]{\ensuremath{\,(\mathrm{mod}\,#1)}}

\newcommand{\Matrix}[1]{\begin{pmatrix}#1\end{pmatrix}}
\newcommand{\SmallMatrix}[1]{\left(\begin{smallmatrix}#1\end{smallmatrix}\right)}

\renewcommand{\Re}{{\mathop{\mathgroup\symoperators Re}}}
\renewcommand{\Im}{{\mathop{\mathgroup\symoperators Im}}}
\newcommand{\sgn}{{\mathop{\mathgroup\symoperators \,sign}}}

\newcommand{\Max}[1]{\ensuremath{\max \set{#1}}}
\newcommand{\Min}[1]{\ensuremath{\min \set{#1}}}

\newcommand{\Z}{\mathbb{Z}}
\newcommand{\R}{\mathbb{R}}
\newcommand{\N}{\mathbb{N}}
\newcommand{\Q}{\mathbb{Q}}
\newcommand{\C}{\mathbb{C}}

\newcommand{\wbar}[1]{\overline{#1}}
\newcommand{\wtilde}[1]{\widetilde{#1}}
\newcommand{\what}[1]{\widehat{#1}}
\newcommand{\BigO}[2][]{O_{#1}\paren{#2}}
\newcommand{\e}[1]{e\paren{#1}}
\newcommand{\trans}[1]{{^t #1}}

\newcommand{\derv}[3][]{{\frac{d^{#1}#2}{d{#3}^{#1}}}}
\newcommand{\pderv}[3][]{{\frac{\partial^{#1}#2}{\partial{#3}^{#1}}}}

\theoremstyle{plain} 
\newtheorem{thm}{Theorem}

\newif\ifshowTODOs
\showTODOsfalse

\newcommand{\cosmu}{\mathbf{cos}}
\newcommand{\sinmu}{\mathbf{sin}}
\newcommand{\specmu}{\mathbf{spec}}

\newcommand{\textdef}[1]{\textit{#1}}
\DeclareMathOperator*{\res}{res}

\title{The Spectral Kuznetsov Formula on $SL(3)$}
\author{Jack Buttcane}
\date{14 November 2014}
\address{Mathematisches Institut, Georg-August Universit\"at G\"ottingen, Bunsenstrasse 3-5, D-37073 G\"ottingen, Germany}
\email{buttcane@uni-math.gwdg.de}
\thanks{During the time of this research, the author was supported by European Research Council Starting Grant agreement number 258713.}

\begin{document}

\begin{abstract}
The $SL(3)$ Kuznetsov formula exists in several versions, and has been employed with some success to study automorphic forms on $SL(3)$.
In each version, the weight functions on the geometric side are given by multiple integrals with complicated oscillating factors; this is the primary obstruction to its use.
By describing them as solutions to systems of differential equations, we give power series and Mellin-Barnes integral representations of minimal dimension for these weight functions.
This completes the role of harmonic analysis on symmetric spaces on the geometric side of the Kuznetsov formula, so that further study may be done through classical analytic techniques.
\end{abstract}

\subjclass[2010]{Primary 11F72; Secondary 44A20, 33E20}

\maketitle

The classical Kuznetsov formulas serve as a bridge between Fourier coefficients of Maass forms and Kloosterman sums.
They are a necessary step in most proofs of subconvexity for $L$-functions, in a variety of equidistribution problems, in proofs of bounds and asymptotics for Maass forms, and in many applications of the circle method to problems in additive number theory.
In attempting to generalize the Kuznetsov formula to study Fourier coefficents of $SL(3)$ Maass forms, we encounter difficulties due to the lack of knowledge of the generalized Bessel functions that appear.
The goal of this paper is to solve those difficulties for the spectral Kuznetsov formula.

There are essentially four versions of the spectral Kuznetsov formula:
The first, due to Xiaoqing Li, appears in \cite[Thm 11.6.19]{Goldfeld} and uses an idea of Don Zagier with bi-$K$ invariant test functions and spherical inversion.
The second, due to Valentin Blomer, appears in \cite{Bl01} and uses the classical construction with the inner product of two Poincar\'e series.
In \cite{GK}, Dorian Goldfeld and Alex Kontorovich applied Lebedev-Whittaker inversion to Blomer's construction to obtain an arbitrary test function on the spectrum.
Finally, in \cite{Me01}, the author worked out the first part of the integral transform on the geometric side of Li's version to obtain a slightly more usable formula -- the main thrust there being an approximation to the geometric Kuznetsov formula.

These four versions have been used with varying degrees of success.
Li was able to obtain some results on spectral averages of the $1,1$ Fourier coefficients of cusp forms, weighted by their Whittaker functions \cite{Li01}.
Blomer has given a non-optimal version of the spectral large sieve inequalities and Lindel\"of on average for the second moment of $SL(3)$ $L$-functions \cite{Bl01}.
Blomer, Raulf, and the author have given results on the distribution of the Fourier coefficients of $SL(3)$ forms, including a Sato-Tate law on average \cite{SatoTateLaw}.
Kontorovich and Goldfeld obtained some results on low-lying zeros and symmetry types \cite{GK}.
In every case, the treatment of the weight functions was a significant challenge in the analysis.

If we have a nice test function $f(\mu)$, $\mu=(\mu_1,\mu_2,\mu_3)\in\C^3$ with $\mu_1+\mu_2+\mu_3=0$, which is invariant under permutations of the coordinates of $\mu$, then the $SL(3)$ Kuznetsov formula is formally
\begin{align*}
	& \sum_{\varphi} f(\mu_\varphi) \rho_\varphi^*(n) \wbar{\rho_\varphi^*(m)} + \text{ Eisenstein series terms} \\
	&= \delta_{m=n} H_I(f;(-1, -1)) + \text{ intermediate Weyl element terms} \\
	& \qquad + \sum_{\varepsilon\in\set{\pm1}^2} \sum_{c_1,c_2\in\N} \frac{S_{w_l}(\psi_m,\psi_{\varepsilon n},c)}{c_1 c_2} H_{w_l}\paren{f; \paren{-\varepsilon_2 \tfrac{m_1 n_2 c_1}{c_2^2}, -\varepsilon_1 \tfrac{m_2 n_1 c_2}{c_1^2}}}, \nonumber
\end{align*}
where $\set{\varphi}$ is a basis of Maass cusp forms with suitably normalized Fourier coefficients $\rho_\varphi^*(n)$ and spectral parameters $\mu_\varphi$, and $S_{w_l}(\psi_m,\psi_n,c)$ is the generalized Kloosterman sum attached to the long Weyl element.
In Li's version of the spectral $SL(3)$ Kuznetsov formula, the weight functions $H_w\paren{f;y}$ are given by a rather complicated integral transform:
Let $h_\mu$ be the $SL(3,\R)$ spherical function, and $y=(y_1,y_2)\in \R^2$ with $y_1 y_2 \ne 0$, then the long element weight function is
\begin{align}
\label{eq:HwIntro}
	H_{w_l}(f; y) =& \frac{1}{48 \pi^4 \abs{y_1 y_2}} \int_0^\infty \int_0^\infty \int_{\R^3} \int_{\R^3} \int_{\Re(\mu)=0} f(\mu) \e{x_1+x_2+x_1'+x_2'} \\
	& \quad h_\mu\paren{\SmallMatrix{t_1 t_2\\&t_2\\&&1}^{-1} \SmallMatrix{1&x_2&x_3\\&1&x_1\\&&1}^{-1} \SmallMatrix{y_1 y_2\\&y_2\\&&1} \SmallMatrix{&&1\\&-1\\1} \SmallMatrix{1&x_2'&x_3'\\&1&x_1'\\&&1} \SmallMatrix{t_1 t_2\\&t_2\\&&1}} \nonumber \\
	& \quad \specmu(\mu) d\mu \, dx \, dx' \frac{dt_1 \, dt_2}{t_1 t_2^2}, \nonumber
\end{align}
where $\e{x} = e^{2\pi i x}$, $d\mu = d\mu_1 \, d\mu_2$, and $\specmu(\mu)$ is the spectral measure
\[ \specmu(\mu) = -\prod_{j<k} (\mu_j-\mu_k) \tan \frac{\pi}{2}(\mu_j-\mu_k). \]
Aside from the large matrix product in the middle, this is a 10-dimensional integral transform -- 13-dimensional if we count the definition of $h_\mu$, and the combined integral does not converge absolutely.

Providing good bounds for this integral transform has been quite difficult.
A singularly vexing quantity which appears in every version of the formula (after some substitutions) is the oscillating factor
\begin{align}
\label{eq:Psixstar}
	\e{t_1 x_1'+t_2 x_2'-\frac{y_1}{t_2}\frac{x_2'+x_1' x_3'}{1+{x_2'}^2+{x_3'}^2}-\frac{y_2}{t_1}\frac{x_1'+x_2'(x_1'x_2'-x_3')}{1+{x_1'}^2+(x_1' x_2'-x_3')^2}}.
\end{align}
Attempting to apply arguments like stationary phase to this factor results in a combinatorial mess, and efforts to evaluate the integrals in \eqref{eq:HwIntro} using classical special functions seem to stall out when treating it.

There are four integral transforms in total, corresponding to the three non-degenerate Weyl elements plus the trivial term, which appear in the Kuznetsov formula.
The general form of the transform is given in \eqref{eq:HwDef} below, of which \eqref{eq:HwIntro} is a particular case.
The four transforms are
\[ H_I(f;(-1,-1)), \qquad H_{w_4}(f;(y_1,-1)), \qquad H_{w_5}(f;(-1,y_2)), \qquad H_{w_l}(f;(y_1,y_2)), \]
and it is the main theorem of this paper that they may be given as integrals of the test function $f$ against kernel functions with simple power series expansions.
We denote the frequently used functions
\[ \cosmu(\mu) = \prod_{j<k} \cos \frac{\pi}{2}(\mu_j-\mu_k), \qquad \sinmu(\mu) = \prod_{j<k} \sin \frac{\pi}{2}(\mu_j-\mu_k). \]
\begin{thm}
	Suppose $f(\mu)$ is symmetric in $\mu$, holomorphic in a neighborhood of $\Re(\mu)=0$, and Schwartz class in each $\Im(\mu_i)$, then for $y_1,y_2\in \R\setminus\set{0}$, we have
\[ H_{w_5}(f; (-1,y_2)) = H_{w_4}(\wtilde{f}; (-y_2,-1)), \qquad \wtilde{f}(\mu) = f(-\mu), \]
\begin{align*}
	H_{w_4}(f; (y_1,-1)) =& \frac{1}{2^{14}\pi^6 \abs{y_1}} \int_{\Re(\mu)=0} f(\mu) \, J_{w_4}(y_1, \mu) \sin\frac{\pi}{2}(\mu_1-\mu_2) \frac{\specmu(\mu)}{\sinmu(\mu)} d\mu,
\end{align*}
\[ H_{w_l}(f; y) = -\frac{1}{2^9 \pi^3 \abs{y_1 y_2}} \int_{\Re(\mu)=0} f(\mu) \, J_{w_l}(y, \mu) \frac{\specmu(\mu)}{\sinmu(\mu)} d\mu, \]
where
\[ J_{w_4}(y_1,\mu) = \abs{8\pi^3 y_1}^{1-\mu_3} \sum_{n=0}^\infty \frac{(8\pi^3 i y_1)^n}{n! \, \Gamma\paren{n+1+\mu_1-\mu_3} \, \Gamma\paren{n+1+\mu_2-\mu_3}}, \]
\[ J_{w_l}(y,\mu) = \abs{4\pi^2 y_1}^{1-\mu_3} \abs{4\pi^2 y_2}^{1+\mu_1} \sum_{n_1,n_2\ge 0} \frac{\Gamma\paren{n_1+n_2+\mu_1-\mu_3+1} \, (4\pi^2 y_1)^{n_1} (4\pi^2 y_2)^{n_2}}{\prod_{i=1}^3 \Gamma\paren{n_1+\mu_i-\mu_3+1}\Gamma\paren{n_2+\mu_1-\mu_i+1}}. \]
\end{thm}
The trivial term was previously evaluated in \cite{Me01}:
\[ H_I(f;(-1,-1)) = \frac{1}{192 \pi^5} \int_{\Re(\mu)=0} f(\mu) \, \specmu(\mu) d\mu. \]
Note that, in practice, the decay assumptions on $f$ can be relaxed to absolute convergence of both sides of the Kuznetsov formula.
We will state the full spectral Kuznetsov formula in section \ref{sect:LiKuz}.

We obtain these results indirectly by showing that the integrals of \eqref{eq:HwIntro} can be rearranged into an integral of the test function $f$ against a kernel function (also defined by a complicated integral), and then showing the kernel function must satisfy a system of differential equations, given explicitly by \eqref{eq:XD1and2}-\eqref{eq:XD2} and \eqref{eq:Kw4Delta3}-\eqref{eq:Kw4Delta3Def}.
Then we find the set of solutions to the system and the particular linear combination that gives the kernel function.
These steps are not difficult from a heuristic vantage point, but justifying them is somewhat delicate.

It is interesting to note that the differential equations we obtain in the long element case are strongly related to those satisfied by the Whittaker function.
To be precise, a function $f_\mu(y_1, y_2)$ satisfies the differential equations of the Whittaker function with spectral parameters $\frac{1}{2}\mu$ if and only if $\sqrt{y_1 y_2} f_\mu(2\sqrt{y_1},2\sqrt{y_2})$ satisfies the differential equations \eqref{eq:XD1and2}-\eqref{eq:XD2}.
This fact will appear in Theorem \ref{thm:LongEleMBs} part \ref{itm:LongElePosCase}, below.
A word of warning concerning the comparison:
Because of the absolute values in the definition of $J_{w_l}$, the function at negative values of the $y_i$ is not the analytic continuation of the function at positive values (they differ by a function of $\mu$).
In particular, the other sign cases are not given by analytic continuation of the Whittaker function.

Because the various integral representations are key in applications, we list the Mellin-Barnes integral forms for the weight functions at each sign.
The function which has a good Mellin-Barnes representation varies based on the signs of each $y_i$.

The Weyl group $W$ acts on $\mu$ by permutations, which we denote $\mu^w$ for $w\in W$.
Let $W_3$ be the subgroup of $W$ generated by the order-three Weyl elements.
Now consider the more symmetric functions
\begin{align*}
	K_{w_l}(y, \mu) =& -\frac{\pi^3}{32} \sum_{w\in W} \frac{1}{\sinmu(\mu^w)} J_{w_l}(y,\mu^w), \\
	K_{w_l}^{+-}(y, \mu) =& J_{w_l}(y,(\mu_1,\mu_2,\mu_3)) - J_{w_l}(y,(\mu_1,\mu_3,\mu_2)), \\
	K_{w_l}^{-+}(y, \mu) =& J_{w_l}(y,(\mu_1,\mu_2,\mu_3))-J_{w_l}(y,(\mu_2,\mu_1,\mu_3)), \\
	K_{w_l}^{--}(y, \mu) =& J_{w_l}(y,(\mu_1,\mu_2,\mu_3)) - J_{w_l}(y,(\mu_3,\mu_2,\mu_1)), \\
	K_{w_4}(y_1,\mu) =& \frac{1}{512} \sum_{w\in W_3} \frac{1}{\sin\frac{\pi}{2}(\mu_1^w-\mu_3^w) \sin\frac{\pi}{2}(\mu_2^w-\mu_3^w)} J_{w_4}(y_1,\mu^w).
\end{align*}
These versions of the kernel functions occur as
\begin{align*}
	H_{w_l}(f; y) =& \frac{1}{96 \pi^6 \abs{y_1 y_2}} \int_{\Re(\mu)=0} f(\mu) \, K_{w_l}(y, \mu) \, \specmu(\mu) d\mu \\
	=& -\frac{1}{2^{10} \pi^3 \abs{y_1 y_2}} \int_{\Re(\mu)=0} f(\mu) \, K_{w_l}^{\pm\pm}(y, \mu) \frac{\specmu(\mu)}{\sinmu(\mu)} d\mu, \\
	H_{w_4}(f; (y_1,-1)) =& \frac{1}{96 \pi^6 \abs{y_1}} \int_{\Re(\mu)=0} f(\mu) \, K_{w_4}(y_1, \mu) \, \specmu(\mu) d\mu.
\end{align*}

\begin{thm}
\label{thm:LongEleMBs}
	For the long element weight function:
\begin{enumerate}
\item\label{itm:LongElePosCase} If $y_1,y_2 > 0$,
\begin{align*}
	K_{w_l}(y,\mu) = \pi^4 \cosmu(\mu) \sqrt{y_1 y_2} \; W^*\paren{(2 \sqrt{y_1},2\sqrt{y_2}),2\mu,\psi_{1,1}},
\end{align*}
where $W^*(y,\mu,\psi_{1,1})$ is the completed Whittaker function, which has the Mellin-Barnes integral representation \eqref{eq:WhittakerMellin}.

\item If $y_1,y_2 < 0$,
\begin{align*}
	K_{w_l}^{--}(y, \mu) &= -\frac{1}{\pi} \sin\pi(\mu_1-\mu_3) \int_{-i\infty}^{+i\infty} \int_{-i\infty}^{+i\infty} \abs{4\pi^2 y_1}^{1-s_1} \abs{4\pi^2 y_2}^{1-s_2} \\
	& \qquad \times \frac{\Gamma\paren{s_1-\mu_3} \Gamma\paren{s_1-\mu_1} \Gamma\paren{s_2+\mu_1} \Gamma\paren{s_2+\mu_3}}{\Gamma\paren{1-s_1+\mu_2} \Gamma\paren{s_1+s_2} \Gamma\paren{1-s_2-\mu_2}} \frac{ds_1}{2\pi i} \frac{ds_2}{2\pi i}.
\end{align*}

\item If $y_1 < 0 < y_2$,
\begin{align*}
	K_{w_l}^{-+}(y, \mu) &= -\frac{1}{\pi} \sin\pi(\mu_1-\mu_2) \int_{-i\infty}^{+i\infty} \int_{-i\infty}^{+i\infty} \abs{4\pi^2 y_1}^{1-s_1} \abs{4\pi^2 y_2}^{1-s_2} \\
	& \qquad \times \frac{\Gamma\paren{s_1-\mu_3} \Gamma\paren{s_2+\mu_1} \Gamma\paren{s_2+\mu_2} \Gamma\paren{1-s_2-s_1}}{\Gamma\paren{1+\mu_1-s_1} \Gamma\paren{1+\mu_2-s_1} \Gamma\paren{1-s_2-\mu_3}} \frac{ds_1}{2\pi i} \frac{ds_2}{2\pi i}.
\end{align*}

\item If $y_1 > 0 > y_2$,
\begin{align*}
	K_{w_l}^{+-}(y, \mu) &= -\frac{1}{\pi} \sin\pi(\mu_2-\mu_3) \int_{-i\infty}^{+i\infty} \int_{-i\infty}^{+i\infty} \abs{4\pi^2 y_1}^{1-s_1} \abs{4\pi^2 y_2}^{1-s_2} \\
	& \qquad \times \frac{\Gamma\paren{s_1-\mu_2} \Gamma\paren{s_1-\mu_3}\Gamma\paren{s_2+\mu_1} \Gamma\paren{1-s_1-s_2}}{\Gamma\paren{1-s_1+\mu_1}\Gamma\paren{1-\mu_2-s_2} \Gamma\paren{1-\mu_3-s_2}} \frac{ds_2}{2\pi i} \frac{ds_1}{2\pi i}.
\end{align*}

\end{enumerate}
\end{thm}
A word on contours:
For integrals of the type in Theorem \ref{thm:LongEleMBs}, we generally follow the Barnes integral convention that the contour should pass to the right of all of the poles of the gamma functions in the form $\Gamma(s_i+a)$ and to the left of all of the poles of the gamma functions in the form $\Gamma(a-s_i)$.
This is sufficient for all but those integrals whose integrands lack a net exponential decay; the integrals in parts (2), (3) and (4) have this difficulty.
By Stirling's approximation for the gamma function, we need the unbounded portion of the contour in both $s_i$ to pass some distance to the left of the imaginary axis for the integrals to converge.
For $\max_i \abs{\Re(\mu_i)}<\eta$ with some small $\eta > 0$, we may set each contour at $\Re(s_i) = \eta$ for $\abs{\Im(s_i)} \le \max_i \abs{\Im(\mu_i)}+\eta$ and $\Re(s_i) = -2\eta$ for $\abs{\Im(s_i)} \ge \max_i \abs{\Im(\mu_i)}+\eta$, connected by horizontal lines.
This handles $\Re(\mu)$ near zero, otherwise it is likely more useful to simply shift both contours to the vertical lines at, say $\Re(s_1) =\Re(s_2) = -\eta < 0$, picking up the residues at the relevant poles.

\begin{thm}
\label{thm:w4MBs}
	For the $w_4$ weight function:
\begin{align*}
	K_{w_4}(y,\mu) =& \frac{1}{512 \pi^2} \int_{-i\infty}^{+i\infty} \abs{8\pi^3 y_1}^{1-s} \Gamma(s-\mu_1)\Gamma(s-\mu_2)\Gamma(s-\mu_3) \\
	& \qquad \paren{\exp\paren{-i\tfrac{3\pi}{2}\varepsilon s}+\exp\paren{i\tfrac{\pi}{2} \varepsilon s}\sum_{j=1}^3\exp\paren{i\pi\varepsilon \mu_j}} \frac{ds}{2\pi i},
\end{align*}
	where $\varepsilon = \sgn(y_1)$.
\end{thm}
The integral of Theorem \ref{thm:w4MBs} also lacks a net exponential decay, but has sufficient polynomial decay if we take, say $\max_i \Re(\mu_i) < \Re(s) < \frac{1}{6}$.
The proofs of these representations are given in section \ref{sect:MBRepns}.

We note that the integral representation of Theorem \ref{thm:w4MBs} matches that of the function $g^\pm(y)$ in the proof of \cite[Lemma 6]{Bl02}.
(The formula for $G^\pm(s)$ there should contain the sum of the three terms in the innermost parentheses, as opposed to the product.)
This gives another Mellin-Barnes representation of $K_{w_4}$,
\begin{align}
\label{eq:Kw4Voronoi}
	K_{w_4}(y,\mu) =& \frac{1}{128 \sqrt{\pi}} \int_{-i\infty}^{+i\infty} \abs{\pi^3 y_1}^{1-s} \paren{\prod_{j=1}^3 \frac{\Gamma\paren{\frac{s-\mu_j}{2}}}{\Gamma\paren{\frac{1-s+\mu_j}{2}}}+\varepsilon i \prod_{j=1}^3 \frac{\Gamma\paren{\frac{1+s-\mu_j}{2}}}{\Gamma\paren{\frac{2-s+\mu_j}{2}}}} \frac{ds}{2\pi i}.
\end{align}
This is the weight function occuring in the Voronoi summation formula for $GL(3)$, which is also attached to the $w_5$ element in some sense (and hence the $w_4$ element, as well).

\section{Background}
\subsection{Groups, spaces and characters}
Let $G=SL(3,\R)$, $\Gamma=SL(3,\Z)$, and $K=SO(3,\R)$.
The Weyl group $W$ of $G$ contains the six matrices
\begin{equation*}
	\begin{array}{rclcrclcrcl}
		I &=& \Matrix{1\\&1\\&&1}, && w_2 &=& \Matrix{&1\\-1\\&&1}, && w_3 &=& \Matrix{1\\&&-1\\&1}, \\
		w_4 &=& \Matrix{&1\\&&1\\1}, && w_5 &=& \Matrix{&&1\\1\\&1}, && w_l &=& \Matrix{&&1\\&-1\\1},
	\end{array}
\end{equation*}
where the missing entries should be interpreted as zero.
We will generally not be concerned with the elements $w_2$ and $w_3$, as they don't appear in the Kuznetsov formula for non-degenerate characters.

We denote the spectral parameters of Maass forms as $\mu = (\mu_1,\mu_2,\mu_3) \in \C^3$ subject to $\sum_i \mu_i = 0$.
An integral over the space of such $\mu$, e.g. at $\Re(\mu)=0$, will be denoted
\[ \int_{\Re(\mu)=0} f(\mu) d\mu := \int_{(0)} \int_{(0)} f(\mu_1,\mu_2,-\mu_1-\mu_2) d\mu_1 \, d\mu_2. \]
Because the Jacobian for the change of variables has determinant $\pm1$ in every case, we could have chosen any of
\[ d\mu := d\mu_1 \, d\mu_2 = d\mu_1 \, d\mu_3 = d\mu_2 \, d\mu_3. \]
In a similar manner, we tend to denote multivariable inverse Mellin transforms as
\[ \int_{\Re(s) = (c_1, \ldots, c_n)} f(s) \frac{ds}{(2\pi i)^n} := \int_{(c_1)} \cdots \int_{(c_n)} f(s_1, \ldots, s_n) \frac{ds_1 \cdots ds_n}{(2\pi i)^n}. \]
The exception to this will be when writing an integral of the form
\[ \int_{-i\infty}^{+i\infty} f(s) \frac{ds}{2\pi i}. \]
In this case, $f(s)$ will involve a quotient of gamma functions and we follow the Barnes integral convention.

The space of $x$-matrices for a ring $R \in \set{\R,\Q,\Z}$ is denoted
\[ U(R) = \set{\Matrix{1&x_2&x_3\\&1&x_1\\&&1} \setdiv x_i\in R}. \]
The relationship between the indices of the $x_i$ and their location within the $x$-matrix will be fixed throughout the paper.
The measure on the space $U(\R)$ is simply $dx = dx_1 \, dx_2 \, dx_3$.
We will tend to denote $x$-matrices by the symbols $x$, $x'$, $x^*$ or $u$, except as they pertain to Kloosterman sums, where we use $b$ or $b'$.
Characters of $U(\R)$ are given by
\[ \psi_m(x) = \psi_{m_1,m_2}(x) = \e{m_1 x_1+m_2 x_2}, \]
where $m\in\R^2$; the symbol $\psi$ will generally denote such a character.
We will tend to assume $m \in \Z^2$, so that $\psi_m$ is trivial on $U(\Z)$

The space $U$ decomposes into subspaces according to the action of the Weyl group, and we set $\wbar{U}_w = (w^{-1} \, \trans{U} \, w) \cap U$ for $w\in W$.
For the particular Weyl elements $I, w_4, w_5$ and $w_l$, these are
\begin{align}
\label{eq:ExplicitUw}
\wbar{U}_I = \set{I}, \qquad \wbar{U}_{w_4} = \set{\Matrix{1&x_2&x_3\\&1&0\\&&1}}, \qquad \wbar{U}_{w_5} = \set{\Matrix{1&0&x_3\\&1&x_1\\&&1}}, \qquad \wbar{U}_{w_l} = U.
\end{align}
The complementary spaces $U_w = (w^{-1} \, U \, w) \cap U$ will make a brief appearance in section \ref{sect:KloosSums}.

Define the spaces of $y$-matrices as
\[ Y = \set{\Matrix{y_1 y_2\\& y_1\\&&1} \in GL(3,\R)}, \qquad Y^+ = \set{y \in Y \setdiv y_1,y_2 > 0}, \]
with measure
\[ dy := \frac{dy_1 \, dy_2}{(y_1 y_2)^3}. \]
We tend to use the symbols $y$, $y^*$ or $t$ to denote $y$-matrices, and we do not distinguish notationally between the matrix $\SmallMatrix{y_1 y_2\\& y_1\\&&1}$ and the pair $(y_1, y_2)$.
This will not cause a problem, as the multiplication is the same in both realizations.

The Weyl group acts on $y$-matrices by
\[ y^w = w y w^{-1} \pmod{\pm \R^+}; \]
in particular, we have
\begin{align}
\label{eq:ExplicitWeylActy}
	y^I = y, \qquad y^{w_4} = \paren{\frac{1}{y_1 y_2},y_1}, \qquad y^{w_5} = \paren{y_2, \frac{1}{y_1 y_2}}, \qquad y^{w_l} = \paren{\frac{1}{y_2}, \frac{1}{y_1}}.
\end{align}

Now a matrix $g$ in $GL(3,\R)$ has the \textdef{Iwasawa decomposition} $g=rxyk$ with $r \in \R\setminus\set{0}$, $x\in U(\R)$, $y\in Y$ and $k \in K$.
In this decomposition, we have a choice of signs for the $y_i$.
When applying the decomposition to $g\in G$, unless otherwise noted, we will assume that $y \in Y^+$ and hence $r\in \R^+$.

On the space of diagonal matrices in $GL(3,\R)$, for $\mu$ as above, we define the power function
\[ p_\mu\paren{\Matrix{a_1\\&a_2\\&&a_3}} = \prod_{i=1}^3 \abs{a_i}^{\mu_i}, \]
and extend to all of $GL(3,\R)$ by the Iwasawa decomposition
\[ p_\mu(rxyk) = p_\mu(y). \]
Note that this is well-defined over the choice of signs for the $y_i$.
Set $\rho=(1,0,-1)$, then we will frequently consider the normalized power function
\[ p_{\rho+\mu}(rxyk) = \abs{y_1}^{1-\mu_3} \abs{y_2}^{1+\mu_1}. \]
The Weyl group acts as permutations on the coordinates of $\mu$ via the power function
\[ p_{\mu^w}(y) = p_\mu(y^w) = p_\mu(w y w^{-1}); \]
in particular, we have
\begin{align}
\label{eq:ExplicitWeylActmu}
	\mu^I = \mu, \qquad \mu^{w_4} = (\mu_3,\mu_1,\mu_2), \qquad \mu^{w_5} = (\mu_2,\mu_3,\mu_1), \qquad \mu^{w_l} = (\mu_3,\mu_2,\mu_1).
\end{align}
We say that a function is \textdef{symmetric in $\mu$} when it is invariant under permutations of the coordinates of $\mu$.

\subsection{The Lie algebra of $G$}
We need a very small amount of information on Lie algebras:
Denote by $E_{ij}$ the $3 \times 3$ matrix with a one in the $i$-th row and $j$-th column, and zeros elsewhere, then these matrices form a basis of the Lie algebra of $G$, and they act on smooth functions on $G$ by
\begin{align}
\label{eq:LieAlgebraAction}
	E_{ij}f(g) =& \left.\pderv{}{t}f(g (I+t E_{ij}))\right|_{t=0}.
\end{align}

A basis for the center of the Lie algebra is given by
\begin{align}
\label{eq:Delta1Def}
	\Delta_1 =& \frac{1}{2} \sum_{i,j} E_{ij} E_{ji}, \\
\label{eq:Delta2Def}
	\Delta_2 =& \frac{1}{3} \sum_{i,j,k} E_{ij} E_{jk} E_{ki} - \Delta_1,
\end{align}
and, acting on functions in $C^\infty(G)$, these operators are left- and right-translation invariant.
This last is a key point in what follows, and is not true if we only consider the $\Delta_i$ as operators on $C^\infty(G/K)$.
In particular, we will not be able to use the representations of the operators given in \cite[eq. (6.1.1)]{Goldfeld}.

\subsection{Bessel functions}
We may write the $K$-Bessel function as the Mellin-Barnes integral \cite[6.561.16]{GR}
\begin{align}
\label{eq:KBesselMellin}
	K_\nu(2x) =& \frac{1}{4} \int_{\Re(s)=1} \Gamma\paren{\frac{s+\nu}{2}} \Gamma\paren{\frac{s-\nu}{2}} x^{-s} \frac{ds}{2\pi i},
\end{align}
for $\abs{\Re(\nu)} < 1$.

\subsection{Whittaker functions}
The incomplete Whittaker function is defined by the Jacquet integral formula
\begin{align}
\label{eq:JacquetWhitt}
	W(z,\mu,\psi_m) =& \int_{U(\R)} p_{\rho+\mu}(w_l uz) \psi_m(u) du.
\end{align}
Its completion is given by $W^*(z;\mu,\psi_m) = \Lambda(\mu) W(z;\mu,\psi_m)$, where
\[ \Lambda(\mu) = \pi^{-\frac{3}{2}+\mu_3-\mu_1} \Gamma\paren{\frac{1+\mu_1-\mu_2}{2}} \Gamma\paren{\frac{1+\mu_1-\mu_3}{2}} \Gamma\paren{\frac{1+\mu_2-\mu_3}{2}}, \]
and is symmetric in $\mu$.
We will use several times the fact that
\begin{align}
\label{eq:LambdaCos}
	\Lambda(\mu)\Lambda(-\mu)\cosmu(\mu) = 1,
\end{align}
which follows from the reflection formula \eqref{eq:ReflectionDuplication}.

The Mellin transform of the completed Whittaker function is known \cite[eqs. (6.1.4), (6.1.5)]{Goldfeld}
\begin{align}
\label{eq:WhittakerMellin}
	W^*(y,\mu,\psi_{1,1}) &= \frac{1}{4\pi^2} \int_{\Re(s) = (2,2)} G(s, \mu) (\pi y_1)^{1-s_1} (\pi y_2)^{1-s_2} \frac{ds}{(2\pi i)^2},
\end{align}
where
\[ G(s, \mu) = \frac{\Gamma\paren{\frac{s_1-\mu_1}{2}}\Gamma\paren{\frac{s_1-\mu_2}{2}}\Gamma\paren{\frac{s_1-\mu_3}{2}}\Gamma\paren{\frac{s_2+\mu_1}{2}}\Gamma\paren{\frac{s_2+\mu_2}{2}}\Gamma\paren{\frac{s_2+\mu_3}{2}}}{\Gamma\paren{\frac{s_1+s_2}{2}}}. \]
From the Mellin transform, assuming the components of $\mu$ are distinct, we may extract the asymptotics
\begin{align}
\label{eq:WhittAsymp1}
	W^*(y,\mu,\psi_{1,1}) \sim& \sum_{w\in W} C_1(\mu^w) p_{\rho+\mu^w}(y), \\
\label{eq:Bdef}
	C_1(\mu) =& \pi^{\mu_1-\mu_3} \prod_{i<j} \Gamma\paren{\frac{\mu_j-\mu_i}{2}},
\end{align}
as $(y_1, y_2) \to 0$, for $\abs{\Re(\mu_i)} < 1$.
The analysis in section \ref{sect:MultConst} will require the first-term asymptotic in $y_1$ alone, which is
\begin{align}
\label{eq:Whitty1Asymp}
	W^*(y,\mu,\psi_{1,1}) \sim& \frac{1}{8\pi^2} \sum_{w\in W_3} \abs{\pi y_1}^{1-\mu_3^w} \abs{\pi y_2}^{1-\frac{\mu_3^w}{2}} \Gamma\paren{\tfrac{\mu_3^w-\mu_1^w}{2}} \Gamma\paren{\tfrac{\mu_3^w-\mu_2^w}{2}} K_{\frac{\mu_1^w-\mu_2^w}{2}}\paren{2\pi \abs{y_2}},
\end{align}
as can be seen from the Mellin-Barnes integral for the $K$-Bessel function \eqref{eq:KBesselMellin}.

The Mellin transform of a product of two Whittaker functions is given by a formula of Stade \cite[Prop 11.6.17]{Goldfeld} (note: the leading 2 should be a $2^{n-1}$ there)
\begin{align}
\label{eq:StadesFormula}
	\int_{Y^+} W^*(t,\mu,\psi_{1,1}) W^*(t,\mu',\psi_{1,1}) \, (t_1^2 t_2)^s \, dt = \frac{1}{4 \pi^{3s} \Gamma\paren{\frac{3s}{2}}} \prod_{j,k} \Gamma\paren{\frac{s+\mu_j+\mu_k'}{2}},
\end{align}
which converges absolutely for $\Re(s) > \max_i \set{\abs{\Re(\mu_i)}+ \abs{\Re(\mu'_i)}}$.

\subsection{Maass forms}
An $SL(3,\Z)$ Maass form $\varphi$ with spectral parameters $\mu_\varphi$ is a smooth, square-integrable function on $\Gamma\backslash G/K$ which shares the eigenvalues of the power function $p_{\rho+\mu_\varphi}$ under $\Delta_1$ and $\Delta_2$.
From the explicit representation of the operators acting on $C^\infty(G/K)$ given in \cite[eq. (6.1.1)]{Goldfeld}, we have
\[ -\Delta_1 p_{\rho+\mu} = \paren{1-\frac{\mu_1^2+\mu_2^2+\mu_3^2}{2}} p_{\rho+\mu},\qquad -\Delta_2 p_{\rho+\mu} = -\mu_1 \mu_2 \mu_3 p_{\rho+\mu}. \]

For a Maass form $\varphi$ with spectral parameters $\mu_\varphi$, we have
\[ \int_{U(\Z)\backslash U(\R)} \varphi(uy) \psi_m(u) du = \frac{\rho_\varphi(m)}{\abs{m_1 m_2}} W^*\paren{\SmallMatrix{\abs{m_1 m_2}\\&\abs{m_1}\\&&1}y,\mu_\varphi,\psi_{1,1}}, \]
for $m_1,m_2\in\Z\setminus\set{0}$, where $\rho_\varphi(m)$ is some constant depending on $\varphi$ and $m$, which we call the \textdef{Fourier-Whittaker coefficient at $m$}.

A Maass form $\varphi$ is a \textdef{cusp form} if
\[ \int_{\wbar{U}_{w_5}(\Z)\backslash\wbar{U}_{w_5}(\R)} \varphi(ug) du = \int_{\wbar{U}_{w_4}(\Z)\backslash\wbar{U}_{w_4}(\R)} \varphi(ug) du = 0, \]
for all $g\in G$.

\subsection{Eisenstein series}
There are two types of Eisenstein series over $SL(3,\Z)$:
\[ E(z;\mu) = \sum_{\gamma \in U(\Z)\backslash \Gamma} p_{\rho+\mu}(\gamma z), \qquad E_\phi(z;\mu_1) = \sum_{\gamma \in P_{2,1} \backslash \Gamma} p_{\mu_1,\phi}(\gamma z), \]
where
\[ P_{2,1} = \set{\Matrix{*&*&*\\ *&*&*\\0&0&1} \in \Gamma}, \]
and $p_{\mu_1,\phi}(rxyk) = (y_1^2 y_2)^{\frac{1}{2}+\mu_1} \phi(x_2+iy_2)$ with $\phi$ any $SL(2,\Z)$ Maass cusp form.
These are $SL(3,\Z)$ Maass forms with spectral parameters $\mu$ and $\paren{\mu_1+\mu_\phi,\mu_1-\mu_\phi,-2\mu_1}$, respectively, where $\frac{1}{4}-\mu_\phi^2$ is the eigenvalue of $\phi$ under the $SL(2,\R)$ Laplacian.

We will denote the Fourier-Whittaker coefficient at $m$ of the Eisenstein series $E(z;\mu)$ by $\rho(m;\mu)$, and similarly $\rho_\phi(m;\mu_1)$ for the Eisenstein series $E_\phi(z;\mu_1)$.

\subsection{Kloosterman sums}
\label{sect:KloosSums}
Let $V$ be the group of diagonal orthogonal matrices in $\Gamma$, then the Bruhat decomposition of some $\gamma\in \Gamma$ takes the form $\gamma=bcvwb'$ with $w\in W$, $v\in V$, $b,b'\in U(\R)$ and $c=\SmallMatrix{\frac{1}{c_2}\\&\frac{c_2}{c_1}\\&&c_1}$ for some $c_1,c_2\in\N$.
The $SL(3,\Z)$ Kloosterman sum attached to the Weyl element $w\in W$ is defined as
\begin{align*}
S_w(\psi_m,\psi_n; c) = \sum_{\gamma = bcwb' \in U(\Z)\backslash \Gamma /V\wbar{U}_w(\Z)} \psi_m(b) \psi_n(b'),
\end{align*}
where $m,n\in\Z^2$.
The quotient by $V$ simply allows us to restrict to positive moduli $c_1$ and $c_2$ by conjugating the $v$ matrix, which contains the signs of the moduli, to the right.

Now the Bruhat decomposition is only defined up to an element of $U_w(\R)$ so we define the Kloosterman sum to be zero unless $\psi_n((cw) u (cw)^{-1}) \psi_m(u^{-1}) = 1$ for all $u\in U_w(\R)$, which we call the compatibility condition.
By a computation of Friedberg \cite[pp.\ 173-174]{Friedberg}, only sums satisfying the compatibility condition occur in the Fourier expansion of a Poincar\'e series.
In particular, only the $I, w_4, w_5$, and $w_l$ elements occur in the Kuznetsov formula for non-degenerate characters.

In \cite{BFG01}, these are given concrete expressions using the two exponential sums
\begin{align*}
	\tilde{S}(m_1,n_1,n_2;D_1,D_2) := \sum_{\substack{C_1 (\text{mod }D_1), C_2 (\text{mod }D_2)\\(C_1,D_1)=(C_2,D_2/D_1)=1}} e\left(n_1\frac{\bar{C_1}C_2}{D_1}+n_2\frac{\bar{C_2}}{D_2/D_1}+m_1\frac{C_1}{D_1}\right),
\end{align*}
for $D_1 | D_2$, and
\begin{align*}
	&S(m_1, m_2, n_1, n_2; D_1, D_2) \\
	& = \sum_{\substack{B_1, C_1 \summod{D_1}\\B_2, C_2 \summod{D_2}\\ }} \e{\frac{m_1B_1 + n_1(Y_1 D_2 - Z_1 B_2)}{D_1} + \frac{m_2B_2 + n_2(Y_2 D_1 - Z_2B_1)}{D_2}},
\end{align*}
where the sum is restricted to
\[ D_1C_2 + B_1B_2 + D_2C_1 \equiv 0 \summod{D_1D_2}, \qquad (B_1, C_1, D_1) = (B_2, C_2, D_2) = 1, \]
and the $Y_i$ and $Z_i$ are defined by
\[ Y_1B_1 + Z_1C_1 \equiv 1 \pmod{D_1}, \qquad Y_2B_2 + Z_2C_2 \equiv 1 \pmod{D_2}. \]
We have
\begin{align*}
	S_{w_6}(\psi_m,\psi_n;c) =& S(n_2,n_1,m_1,m_2;c_1, c_2), \\
	S_{w_5}(\psi_m,\psi_n;c) =& \delta_{\substack{n_1 c_2=m_2 c_1^2 \\ c_1|c_2}} \tilde{S}(n_1,m_1,m_2;c_1,c_2) \\
	S_{w_4}(\psi_m,\psi_n;c) =& \delta_{\substack{n_2 c_1=m_1 c_2^2 \\ c_2|c_1}} \tilde{S}(-n_2,m_2,m_1;c_2,c_1).
\end{align*}
Here $\delta_P$ is one if $P$ is true and zero otherwise.

\subsection{Spherical inversion}
The spherical function over $K$ is the bi-$K$ invariant function on $G$ given by the integral
\[ h_\mu(g) = \int_K p_{\rho+\mu}(k g) dk, \qquad \int_K dk=1. \]
For a smooth, complex-valued function $f(\mu)$ on $\Re(\mu)=0$ of rapid decay in each $\abs{\Re(\mu_i)}$, which is symmetric in $\mu$, we define $F$ on $G$ by the spherical inversion formula of Helgason, Harish-Chandra, and Bhanu-Murthy \cite[sect. 4.3, Thm. 1 and eq. 3.23]{T02}
\[ F(g) = \frac{1}{96 \pi^3} \int_{\Re(\mu)=0} f(\mu) \, h_\mu\paren{g} \, \specmu(\mu) d\mu. \]

\subsection{Li's Kuznetsov formula}
\label{sect:LiKuz}
For $y\in Y$, let $H_w(f; y)$ be given by the iterated integral
\begin{align}
\label{eq:HwDef}
	& \frac{2}{\pi \abs{y_1 y_2}} \int_{Y^+} \int_{\wbar{U}_w(\R)} \int_{U(\R)} F\paren{t^{-1} x^{-1} y w x't} \,\psi_{1,1}(x) \psi_{1,1}(x') dx\,dx'\,t_1^2 t_2\,dt,
\end{align}
where $F(g)$ is the spherical inverse defined above.
Note:  We have changed the $\wbar{\psi_{1,1}(x')}$ in \cite{Me01} to $\psi_{1,1}(x')$ here, which is why we have negatives in the Kuznetsov formula.

We start with Li's version of the $GL(3)$ Kuznetsov formula, after applying spherical inversion and some rearrangement,
\begin{thm}[Li]
\label{thm:LisKuznetsovFormula}
	Let $\set{\varphi}$ be an orthonormal basis of the $SL(3,\Z)$ cusp forms with Langlands parameters $\mu_\varphi$, and $\set{\phi}$ an orthonormal basis of $SL(2,\Z)$ cusp forms with Langlands parameters $\mu_\phi$.
	Let $f(\mu)$ be symmetric in $\mu$, holomorphic on $\abs{\Re(\mu_i)} < \frac{1}{2}+\eta$ for some $\eta > 0$, and Schwartz class in each $\Im(\mu_i)$.
	Then for $m, n \in \N^2$
	\begin{align*}
		& \sum_{\varphi} \frac{f(\mu_\varphi)}{\cosmu(\mu_\varphi)} \rho_\varphi(n) \wbar{\rho_\varphi(m)} \\
		& \qquad +\frac{1}{2\pi i} \sum_{\phi} \int_{\Re(\mu_1)=0} \frac{f\paren{\mu_1+\mu_\phi,\mu_1-\mu_\phi,-2\mu_1}}{\cosmu\paren{\mu_1+\mu_\phi,\mu_1-\mu_\phi,-2\mu_1}} \rho_\phi(n; \mu_1) \wbar{\rho_\phi(m; \mu_1)} d\mu_1 \nonumber \\
		& \qquad +\frac{1}{24 (2\pi i)^2} \int_{\Re(\mu)=0} \frac{f(\mu)}{\cosmu(\mu)} \rho(n; \mu) \wbar{\rho(m; \mu)} d\mu \nonumber \\
		&= \delta_{m=n} H_I(f;(-1, -1)) \nonumber \\
		& \qquad + \sum_{\varepsilon\in\set{\pm1}} \sum_{\substack{c_1,c_2\in\N\\ m_2 c_1=n_1c_2^2}} \frac{S_{w_4}(\psi_m,\psi_{(n_1,\varepsilon n_2)},c)}{c_1 c_2} H_{w_4}\paren{f; \paren{\varepsilon \tfrac{m_1 m_2^2n_2}{c_2^3n_1},-1}}, \nonumber \\
		& \qquad + \sum_{\varepsilon\in\set{\pm1}} \sum_{\substack{c_1,c_2\in\N\\m_1c_2=n_2c_1^2}} \frac{S_{w_5}(\psi_m,\psi_{(\varepsilon n_1,n_2)},c)}{c_1 c_2} H_{w_5}\paren{f; \paren{-1,\varepsilon \tfrac{m_1^2 m_2 n_1}{c_1^3 n_2}}}, \nonumber \\
		& \qquad + \sum_{\varepsilon\in\set{\pm1}^2} \sum_{c_1,c_2\in\N} \frac{S_{w_l}(\psi_m,\psi_{\varepsilon n},c)}{c_1 c_2} H_{w_l}\paren{f; \paren{-\varepsilon_2 \tfrac{m_1 n_2 c_1}{c_2^2}, -\varepsilon_1 \tfrac{m_2 n_1 c_2}{c_1^2}}}. \nonumber
	\end{align*}
\end{thm}
Remarks:
\begin{enumerate}
\item Since all Maass forms for $SL(3,\Z)$ are necessarily even, the assumption that $n_i, m_i > 0$ is purely for notational simplicity.
\item The constant on the minimal Eisenstein series above, after accounting for the slight difference of definitions, does not match that of \cite[Thm 10.13.1]{Goldfeld}; we have included some justification for our constant in section \ref{sect:EisenConsts}.
\end{enumerate}

\subsection{Theorems on gamma functions}
The gamma function satisfies the well-known reflection and duplication formulas
\begin{align}
\label{eq:ReflectionDuplication}
	\Gamma(1-z)\Gamma(z)=\frac{\pi}{\sin(\pi z)}, \qquad \Gamma(z) = \frac{2^{z-1}}{\sqrt{\pi}} \Gamma\paren{\frac{z}{2}} \Gamma\paren{\frac{z+1}{2}}.
\end{align}
In section \ref{sect:MultConst} we will use Barnes' integral lemmas:
\begin{thm}[Barnes' first lemma, {\cite[sect. 1.7]{Bailey}}]
\label{thm:BarnesFirst}
For $a,b,c,d\in\C$,
\begin{align*}
	\frac{1}{2\pi i} \int_{-i\infty}^{+i\infty} \Gamma(a+s) \Gamma(b+s) \Gamma(c-s) \Gamma(d-s) ds &= \frac{\Gamma(a+c) \Gamma(b+c) \Gamma(a+d) \Gamma(b+d)}{\Gamma(a+b+c+d)}.
\end{align*}
\end{thm}
\begin{thm}[Barnes' second lemma, {\cite[sect. 6.2]{Bailey}}]
\label{thm:BarnesSecond}
For $a,b,c,d,e,f\in\C$ with $a+b+c+d+e-f=0$,
\begin{align*}
	& \frac{1}{2\pi i} \int_{-i\infty}^{+i\infty} \frac{\Gamma(a+s) \Gamma(b+s) \Gamma(c+s) \Gamma(d-s) \Gamma(e-s)}{\Gamma(f+s)} ds \\
	&= \frac{\Gamma(a+d) \Gamma(b+d) \Gamma(c+d) \Gamma(a+e) \Gamma(b+e) \Gamma(c+e)}{\Gamma(f-a) \Gamma(f-b) \Gamma(f-c)}.
\end{align*}
\end{thm}

We will also need two Mellin transforms related to the gamma and beta functions:
\begin{align}
\label{eq:BetaDef}
	\int_0^\infty (1+x^2)^u x^{s-1} dx &= \frac{1}{2} B\paren{\frac{s}{2},\frac{-2u-s}{2}},
\end{align}
for $0 < \Re(s) < -2\Re(u)$, which follows from the definition of the beta function $B(u,v) = \frac{\Gamma(u)\Gamma(v)}{\Gamma(u+v)}$, and
\begin{align}
\label{eq:CosMellin}
	\int_0^\infty \cos(x) x^{s-1} dx &= \Gamma(s) \cos\paren{\frac{\pi}{2}s},
\end{align}
for $0<\Re(s)<1$, which follows from the definition of the gamma function.
Note that this last integral does not converge absolutely, but that may be addressed by a simple integration by parts.

\section{The weight functions}
\subsection{The general term}
For the first few steps, we will follow \cite{Me01}, but we simplify the notation somewhat.
In \cite[Lemma 3]{Me01}, the (conditionally convergent) Fourier transform of the spherical function was evaluated
\begin{align}
\label{eq:SphericalFourier}
	\int_{U(\R)} h_\mu(y^{-1} xy') \psi(x) dx = \frac{1}{2\pi^2} W(y,-\mu, \psi) W(y', \mu, \psi),
\end{align}
and also the interchange of $x$ and $\mu$ integrals was justified, so that if
\[ x^* y^* \equiv w x' t \pmod{R^+ K}, \]
then
\begin{align}
\label{eq:GeneralHwStars}
	H_w(f; y) =& \frac{1}{96 \pi^6 \abs{y_1 y_2}} \int_{Y^+} \int_{\wbar{U}_w(\R)} \int_{\Re(\mu)=0} f(\mu) \, W(t,-\mu, \psi_{1,1}) W(y \, y^*, \mu, \psi_{1,1}) \\
	& \qquad \specmu(\mu) d\mu \,\psi_y(x^*) \psi_{1,1}(x') dx'\,t_1^2 t_2\,dt. \nonumber
\end{align}
We may rewrite this slightly as
\begin{align}
\label{eq:GeneralHwFinal}
	H_w(f; y) =& \frac{1}{96 \pi^6 \abs{y_1 y_2}} \int_{Y^+} \int_{\wbar{U}_w(\R)} \int_{\Re(\mu)=0} f(\mu) \, W(x t,-\mu, \psi_{1,1}) W(y w x t, \mu, \psi_{1,1}) \\
	& \qquad \specmu(\mu) d\mu \, dx \, t_1^2 t_2 \, dt, \nonumber
\end{align}
and this is the form we will use to evaluate the $w=I, w_5$ terms.

\subsection{The trivial term}
The identity element weight function is computed explicitly by Stade's formula \eqref{eq:StadesFormula}; when $\mu'=-\mu$ and $s=1$, this becomes
\begin{align*}
	\int_{Y^+} W^*(t,\mu,\psi_{1,1}) W^*(t,-\mu,\psi_{1,1}) \, t_1^2 t_2 \, dt = \frac{\pi}{2 \cosmu(\mu)}.
\end{align*}
From \eqref{eq:GeneralHwFinal} the trivial term is then
\begin{align*}
	H_I(f) = H_I(f; (-1,-1)) =& \frac{1}{192 \pi^5} \int_{\Re(\mu)=0} f(\mu) \, \specmu(\mu) d\mu. \nonumber
\end{align*}

\subsection{\texorpdfstring{The $\boldsymbol{w_5}$ term}{The w5 term}}
Define the involution $\iota$ on $GL(3,\R)$ by $g^\iota = w_l \trans{g^{-1}} w_l$, then the components of the Iwasawa decomposition satisfy
\[ (rxyk)^\iota = r^{-1} x^\iota y^\iota k^\iota, \quad \SmallMatrix{1&x_2&x_3\\&1&x_1\\&&1}^\iota = \SmallMatrix{1&x_1&x_1 x_2-x_3\\&1&x_2\\&&1}, \]
\[ \SmallMatrix{y_1 y_2\\&y_1\\&&1}^\iota \equiv \SmallMatrix{y_2 y_1\\&y_2\\&&1} \pmod{\pm\R^+}, \]
and hence the power function satisfies
\[ p_{\rho+\mu}(g^\iota) = p_{\rho-\mu^{w_l}}(g), \]
and for the Whittaker function
\[ W(g^\iota,\mu, \psi_{1,1}) = W(g,-\mu^{w_l}, \psi_{1,1}). \]
Also the space for the $x'$ integral satisfies $\wbar{U}_{w_5}(\R)^\iota = \wbar{U}_{w_4}(\R)$, and
\[ w_5^\iota = \SmallMatrix{-1\\&-1\\&&1} w_4. \]
Applying these to \eqref{eq:GeneralHwFinal}, we see the $H_{w_5}$ function is given by
\[ H_{w_5}(f; (-1,y_2)) = H_{w_4}(\wtilde{f}; (-y_2,-1)), \qquad \wtilde{f}(\mu) = f(-\mu), \]
so it suffices to determine $H_{w_4}$.

\subsection{The hard terms}
\label{sect:HardTerms}
From the Mellin-Barnes integral representation \eqref{eq:WhittakerMellin}, we notice that the completed Whittaker function may be written
\begin{align}
\label{eq:WhittakerMellinShifted}
	W^*(y,\mu,\psi_{1,1}) &= \frac{1}{4\pi^2} \int_{\Re(s) = (-\eta,-\eta)} G(s, \mu) (\pi y_1)^{1-s_1} (\pi y_2)^{1-s_2} \frac{ds}{(2\pi i)^2} \\
	& \qquad +\frac{1}{4\pi^2} \sum_{w'\in W_3} \int_{\Re(s_2) = -\eta} \paren{\res_{s_1=\mu_3^{w'}} G(s, \mu)} (\pi y_1)^{1-\mu_3^{w'}} (\pi y_2)^{1-s_2} \frac{ds_2}{2\pi i} \nonumber \\
	& \qquad +\frac{1}{4\pi^2} \sum_{w'\in W_3} \int_{\Re(s_1) = -\eta} \paren{\res_{s_2=-\mu_1^{w'}} G(s, \mu)} (\pi y_1)^{1-s_1} (\pi y_2)^{1+\mu_1^{w'}} \frac{ds_1}{2\pi i} \nonumber \\
	& \qquad +\frac{1}{4\pi^2} \sum_{w'\in W} \paren{\res_{s=(\mu_3^{w'}, -\mu_1^{w'})} G(s, \mu)} (\pi y_1)^{1-\mu_3^{w'}} (\pi y_2)^{1+\mu_1^{w'}}, \nonumber
\end{align}
for some small $\eta > 0$, provided the components of $\mu$ are distinct.
Then for the remaining Weyl elements, we apply \eqref{eq:WhittakerMellinShifted} to the second Whittaker function in \eqref{eq:GeneralHwStars}, and substitute $\mu^{w'} \mapsto \mu$ as appropriate to obtain
\begin{align}
\label{eq:HwUsingJw}
	H_w(f; y) =& \frac{1}{96 \pi^6 \abs{y_1 y_2}} \int_{\Re(\mu)=(\eta,0,-\eta)} f(\mu) \, K_w^\text{as}(y, \mu) \, \specmu(\mu) d\mu,
\end{align}
where
\begin{align}
\label{eq:GeneralJw}
	K_w^\text{as}(y, \mu) =& \int_{\Re(s) = (-\eta,-\eta)} \pi^{2-s_1-s_2} G(s, \mu) \, T_w(y,\mu,(-s_2,s_2-s_1,s_1)) \frac{ds}{(2\pi i)^2} \\
	& \quad +3 \int_{\Re(s_2) = -\eta} \pi^{2-\mu_3-s_2} \paren{\res_{s_1=\mu_3} G(s, \mu)} \, T_w(y,\mu,(-s_2,s_2-\mu_3,\mu_3)) \frac{ds_2}{2\pi i} \nonumber \\
	& \quad +3 \int_{\Re(s_1) = -\eta} \pi^{2-s_1+\mu_1} \paren{\res_{s_2=-\mu_1} G(s, \mu)} \, T_w(y,\mu,(\mu_1,-\mu_1-s_1,s_1)) \frac{ds_1}{2\pi i} \nonumber \\
	& \quad +6 \pi^{2-\mu_3+\mu_1} \paren{\res_{s=(\mu_3, -\mu_1)} G(s, \mu)} T_w(y,\mu,\mu), \nonumber
\end{align}
with
\begin{align}
\label{eq:TwDef}
	T_w(y, \mu, \mu') =& \frac{\cosmu(\mu)}{4\pi^2} \int_{Y^+} \int_{\wbar{U}_w(\R)} W^*(t,-\mu, \psi_{1,1}) p_{\rho+\mu'}(y \, y^*) \, \psi_y(x^*) \psi_{1,1}(x') dx'\,t_1^2 t_2\,dt.
\end{align}
Here we have applied \eqref{eq:LambdaCos} and shifted the $\mu$ integrals away from $\Re(\mu)=0$.
The interchange of integrals is justified by absolute convergence of the $T_w$ integral on $\Re(\mu_1) > \Re(\mu_2) > \Re(\mu_3)$, which we will demonstrate for the long element in section \ref{sect:LongEleJustify}.

The notation ``as'' refers to the asymmetry in $\mu$; this is necessary for the moment because the integral defining $T_w$ converges absolutely only on $\Re(\mu_1) > \Re(\mu_2) > \Re(\mu_3)$.
Now $T_w$ extends holomorphically to a neighborhood of $\Re(\mu)=0$, and hence we may discuss the function
\[ K_w(y, \mu) = \frac{1}{6} \sum_{w' \in W} K_w^\text{as}\paren{y, \mu^{w'}}. \]
Again, we will prove this analytic continuation for the long element in section \ref{sect:LongEleJustify}.

We extend $H_w(y,\mu)$ to a function $H_w(g,\mu)$ of $g\in G$, and the construction above allows us to similarly extend $K_w(y, \mu)$ to a function $K_w(g,\mu)$.
Comparing with \eqref{eq:GeneralHwFinal}, we see that $K_w(g, \mu)$ is, in some sense, given by the integral
\begin{align}
\label{eq:JwFakeDef}
	K_w(g, \mu) ``="& \int_{Y^+} \int_{\wbar{U}_w(\R)} W(x t,-\mu, \psi_{1,1}) W(g w x t, \mu, \psi_{1,1}) \, dx \, t_1^2 t_2 \, dt,
\end{align}
though the right-hand side does not converge absolutely for any value of $\mu$.
In sections \ref{sect:LongEleTerm} and \ref{sect:w4Term}, we will treat $K_w(g, \mu)$ as being defined by the integral \eqref{eq:JwFakeDef}, without regard for issues of convergence, and for the long element, we will carefully justify this analysis in section \ref{sect:LongEleJustify}.
Justification for the $w_4$ element analysis is entirely similar.

We wish to point out that our work in section \ref{sect:LongEleJustify} can be used to attach meaning to the integral \eqref{eq:JwFakeDef} directly, but this will not make our analysis any shorter.

The goal for sections \ref{sect:LongEleTerm} and \ref{sect:w4Term} is to describe how the functions $K_w(g, \mu)$ satisfy systems of differential equations, and by solving these equations, determine the $K_w(g, \mu)$ explicitly.

\subsection{The long element term}
\label{sect:LongEleTerm}
As mentioned above, we take \eqref{eq:JwFakeDef} as the definition of $K_w(g,\mu)$.
For the particular case $w=w_l$, this may be written as
\begin{align}
\label{eq:JwlGKDef}
	K_{w_l}(g, \mu) =& \cosmu(\mu) \int_{G/K} W^*(z,-\mu,\psi_{1,1}) W^*(g w_l z,\mu,\psi_{1,1}) \, p_{(1,1,-2)}(z) \, dz.
\end{align}
Again, it is vital to note that the $K_{w_l}$ integral above does not converge absolutely.

\subsubsection{The differential equations}
\label{sect:LongEleDiffEqs}
We wish to analyze $K_{w_l}$ as a function of $g$ via the differential equations it satisfies.
We should point out that this part of the analysis could well have been done before applying \eqref{eq:SphericalFourier}, but for the justification in the following section, our chosen starting point is significantly easier.

The Whittaker function is an eigenfunction of $\Delta_1$ and $\Delta_2$ given by \eqref{eq:Delta1Def} and \eqref{eq:Delta2Def}, with eigenvalues
\begin{align*}
	-\Delta_1 W(g, \mu,\psi_{1,1}) =& \lambda_1 W(g, \mu,\psi_{1,1}), & \lambda_1 = \lambda_1(\mu) =& 1-\tfrac{\mu_1^2+\mu_2^2+\mu_3^2}{2}, \\
	-\Delta_2 W(g, \mu,\psi_{1,1}) =& \lambda_2 W(g, \mu,\psi_{1,1}), & \lambda_2 = \lambda_2(\mu) =& -\mu_1 \mu_2 \mu_3,
\end{align*}
and by right-invariance of the $\Delta_i$, we know that $K_{w_l}(g,\mu)$ also satisfies these equations.

Now for $c,ce-bf \ne 0$ we have the decomposition
\begin{align}
\label{eq:gDecomp}
	s := \Matrix{g&h&i\\d&e&f\\a&b&c} = \Matrix{1&\frac{ch-bi}{ce-bf}&\frac{i}{c}\\&1&\frac{f}{c}\\&&1}\Matrix{\frac{\det s}{ce-bf}\\&\frac{ce-bf}{c}\\&&c}\Matrix{1\\ \frac{cd-af}{ce-bf}\\ \frac{a}{c}&\frac{b}{c}&1},
\end{align}
so away from the degenerate subspace, we may choose coordinates
\begin{align}
\label{eq:gCoords}
	g=&ruyv := r\Matrix{1&u_2&u_3\\&1&u_1\\&&1}\Matrix{y_1 y_2\\&y_1\\&&1} \Matrix{1\\v_1&1\\v_3&v_2&1},
\end{align}
with $r,y_i\in\R\setminus\set{0}$ and $u_i,v_i\in\R$.
Then for functions which are constant in $r$, we apply the definitions \eqref{eq:Delta1Def}, \eqref{eq:Delta2Def} and \eqref{eq:LieAlgebraAction} -- recomputing the decomposition \eqref{eq:gDecomp} at each step, and the differential operators $\Delta_i$ become
\begin{align}
\label{eq:uyvDelta1}
	\Delta_1 =& y_2^2 \partial_{y_2}^2 - y_1 y_2 \partial_{y_1} \partial_{y_2}+y_1^2 \partial_{y_1}^2+y_1(u_2 v_1+y_2) \partial_{u_3} \partial_{v_3}+u_2 y_1 \partial_{u_3} \partial_{v_2}+y_2 \partial_{u_2} \\
	& \qquad +v_1 y_1 \partial_{u_1} \partial_{v_3}+y_1 \partial_{u_1} \partial_{v_2}, \nonumber \\
\label{eq:uyvDelta2}
	\Delta_2 =& -y_2^2 \partial_{y_2}^2+y_1 y_2^2 \partial_{y_1} \partial_{y_2}^2 +y_1^2 \partial_{y_2}^2-y_1^2 y_2 \partial_{y_1}^2 \partial_{y_2}+y_1(u_2 v_1+y_2) \partial_{u_3} \partial_{v_3}+u_2 y_1 \partial_{u_3} \partial_{v_2} \\
	& \qquad +v_1 y_1 y_2 \partial_{u_3} \partial_{v_1} \partial_{v_3}+y_1 y_2 \partial_{u_3} \partial_{v_1} \partial_{v_2}+y_1 y_2 \partial_{u_3}+y_1 y_2 (y_2-u_2 v_1) \partial_{u_3} \partial_{y_2} \partial_{v_3} \nonumber \\
	& \qquad -u_2 y_1 y_2 \partial_{u_3} \partial_{y_2} \partial_{v_2}-y_1^2 y_2 \partial_{u_3} \partial_{y_1} \partial_{v_3}-y_2 \partial_{u_2} \partial_{v_1}+y_1 y_2 \partial_{u_2} \partial_{y_1} \partial_{v_1} \nonumber \\
	& \qquad+u_2 y_1 y_2 \partial_{u_2} \partial_{u_3} \partial_{v_3}+v_1 y_1 \partial_{u_1} \partial_{v_3}+y_1 \partial_{u_1} \partial_{v_2}-v_1 y_1 y_2 \partial_{u_1} \partial_{y_2} \partial_{v_3} \nonumber \\
	& \qquad-y_1 y_2 \partial_{u_1} \partial_{y_2} \partial_{v_2}+y_1 y_2 \partial_{u_1} \partial_{u_2} \partial_{v_3}. \nonumber
\end{align}

Substituting $z \mapsto (w_l v w_l)^{-1} z$ in \eqref{eq:JwlGKDef}, we find
\[K_{w_l}(ruyv, \mu) = \e{u_1+u_2+v_1+v_2} K_{w_l}(y,\mu), \]
so the partials in the $u$ and $v$ variables become
\[ \partial_{u_1} K_{w_l}(g,\mu)=\partial_{u_2} K_{w_l}(g,\mu) = \partial_{v_1}K_{w_l}(g,\mu)=\partial_{v_2}K_{w_l}(g,\mu)=2\pi i K_{w_l}(g,\mu), \]
and
\[ \partial_{u_3} K_{w_l}(g,\mu)=\partial_{v_3} K_{w_l}(g,\mu)=0. \]
The signs on $\partial_{v_1} K_{w_l}(g,\mu)$ and $\partial_{v_2} K_{w_l}(g,\mu)$ are controlled by the inversion and the minus sign in $w_l$.
Thus on the surface $g=\pm1 I y I$, the differential equations on $K_{w_l}(g,\mu)$ become
\begin{align}
\label{eq:XD1and2}
	0 =& \paren{\wtilde{\Delta_i}+\lambda_i(\mu)} K_{w_l}(y,\mu), \\
\label{eq:XD1}
	\wtilde{\Delta_1} =& y_1^2 \partial_{y_1}^2 +y_2^2 \partial_{y_2}^2-y_1y_2\partial_{y_1}\partial_{y_2}+(2\pi i)^2 (y_1+y_2),  \\
\label{eq:XD2}
	\wtilde{\Delta_2} =& -y_1^2 y_2 \partial_{y_1}^2\partial_{y_2}+y_1 y_2^2 \partial_{y_1}\partial_{y_2}^2+y_1 y_2 (2\pi i)^2 \partial_{y_1}-y_1 y_2 (2\pi i)^2 \partial_{y_2}  \\
	& \qquad +y_1^2\partial_{y_1}^2-y_2^2\partial_{y_2}^2+(2\pi i)^2(y_1-y_2) \nonumber.
\end{align}

We follow Bump's analysis of the Whittaker function \cite[ch. II]{B01}:
Consider the operator
\begin{align}
\label{eq:Delta3op}
	\Delta_3 :=& \paren{\Delta_2+\lambda_2(\mu)}+(1-y_1\partial_{y_1})\paren{\Delta_1+\lambda_1(\mu)}.
\end{align}
We know that $\Delta_3\, K_{w_l}(g,\mu)=0$, and, as above, we have $\wtilde{\Delta_3}\, K_{w_l}(y,\mu)=0$, where
\[ \wtilde{\Delta_3} = \paren{\wtilde{\Delta}_2+\lambda_2(\mu)}+(1-y_1\partial_{y_1})\paren{\wtilde{\Delta}_1+\lambda_1(\mu)}. \]
Some rearranging gives
\begin{align*}
	\partial_{y_2} K_{w_l}(y,\mu) =& \frac{1}{(2\pi i)^2 y_1 y_2} \paren{-y_1^3 \partial{y_1}^3- y_1((2\pi i)^2+\lambda_1) \partial_{y_1} + ((2\pi i)^2+\lambda_1+\lambda_2)} K_{w_l}(y,\mu).
\end{align*}
Using this to eliminate the occurances of $\partial_{y_2}$ from \eqref{eq:XD1} gives
\begin{align*}
	0=&\Bigl(y_1^6 \partial_{y_1}^6+6y_1^5\partial_{y_1}^5 + 2 y_1^4 (3-6\pi^2+\lambda_1) \partial_{y_1}^4-2y_1^3(12\pi^2y_1+\lambda_2)\partial_{y_1}^3 \\
	& \qquad +y_1^2(48\pi^4y_1^2 -12\pi^2y_1\lambda_1+3\lambda_1+3\lambda_2)\partial_{y_1}^2 - 2y_1 (3-6\pi^2+\lambda_1)(\lambda_1+\lambda_2)\partial_{y_1} \\
	& \qquad -64\pi^6y_1^2(y_1+y_2) + 16\pi^4y_1^2\lambda_1+2\lambda_1+3\lambda_1(2+\lambda_2) + \lambda_2(6+\lambda_2) \\
	& \qquad -12\pi^2y_1(\lambda_1+\lambda_2)\Bigr) K_{w_l}(y,\mu).
\end{align*}
The two preceeding equations (plus some regularity assumptions) are sufficient to imply that $K_{w_l}(y,\mu)$ belongs to a space of at most 6 linearly independent solutions.

Now assume that the power series
\[ J_{w_l}(y,\mu) = \abs{4\pi^2 y_1}^{1-\mu_3} \abs{4\pi^2 y_2}^{1+\mu_1} \sum_{n_1,n_2\ge 0} a(n_1,n_2) (4\pi^2 y_1)^{n_1} (4\pi^2 y_2)^{n_2}, \]
solves \eqref{eq:XD1and2} (leading power is necessary to make the $n_1=n_2=0$ term work), then the coefficients must satisfy
\begin{align*}
	0 =& \paren{n_1^2-n_1n_2+n_2^2+\mu_2(n_1-2n_2)-\mu_3(n_1+n_2)} a(n_1,n_2) \\
	& \qquad -a(n_1-1,n_2) -a(n_1,n_2-1), \\
	0 =& \paren{(n_1-\mu_3)n_2^2-n_1^2(n_2+\mu_1) -2\mu_2n_1n_2 -\mu_1(\mu_2-\mu_3)n_1 +\mu_3(\mu_2-\mu_1)n_2} a(n_1,n_2) \\
	& \qquad -(n_1-\mu_3) a(n_1,n_2-1)+(n_2+\mu_1) a(n_1-1,n_2),
\end{align*}
where we define $a(n_1,n_2)=0$ when one of the indices is negative.
Adding $(n_2+\mu_1)$ times the first to the second and subtracting $(n_1-\mu_3)$ times the first from the second gives
\begin{align*}
	0 =& \paren{n_1(n_1+\mu_1-\mu_3)(n_1+\mu_2-\mu_3)} a(n_1,n_2) -(n_1+n_2+\mu_1-\mu_3)a(n_1-1,n_2), \\
	0 =& \paren{n_2(n_2+\mu_1-\mu_2)(n_2+\mu_1-\mu_3)} a(n_1,n_2) -(n_1+n_2+\mu_1-\mu_3)a(n_1,n_2-1).
\end{align*}
Thus a solution is given by
\[ J_{w_l}(y,\mu) = \abs{4\pi^2 y_1}^{1-\mu_3} \abs{4\pi^2 y_2}^{1+\mu_1} \sum_{n_1,n_2\ge 0} \frac{\Gamma\paren{n_1+n_2+\mu_1-\mu_3+1} \, (4\pi^2 y_1)^{n_1} (4\pi^2 y_2)^{n_2}}{\prod_{i=1}^3 \Gamma\paren{n_1+\mu_i-\mu_3+1}\Gamma\paren{n_2+\mu_1-\mu_i+1}}, \]
and all 6 linearly independent solutions are given by $J_{w_l}(y,\mu^w)$ for $w \in W$, provided the components of $\mu$ are distinct.

\subsubsection{The asymptotic analysis}
\label{sect:AsympAnal}
We now fix $K_{w_l}$ as a linear combination of these explicit solutions by finding its first-term asymptotic.
To that end, we again write $w_l z \equiv x^* y^* \pmod{\R^+ K}$, so
\begin{align*}
	K_{w_l}(y, \mu) =& \cosmu(\mu) \int_{G/K} W^*(z,-\mu,\psi_{1,1}) W^*(y y^*,\mu,\psi_{1,1}) \psi_y(x^*) \, p_{(1,1,-2)}(z) \, dz.
\end{align*}
Because our analysis takes place inside an integral over $\mu$, we may assume the components of $\mu$ are distinct, so the known asymptotics \eqref{eq:WhittAsymp1} of the Whittaker function on $\Re(\mu)=0$ give
\begin{align*}
	K_{w_l}(y, \mu) \sim& \cosmu(\mu) \sum_{w\in W} C_1(\mu^w) p_{\rho+\mu^w}(y) \int_{G/K} W^*(z,-\mu,\psi_{1,1}) p_{\rho+\mu^w}(y^*) \psi_y(x^*) \, p_{(1,1,-2)}(z) \, dz. \nonumber
\end{align*}
In the limit $(y_1,y_2)\to 0$, the character $\psi_y(x^*) \to 1$, and the $x$ part of the $z$ integral becomes the Jacquet integral \eqref{eq:JacquetWhitt} for the Whittaker function, so
\begin{align*}
	K_{w_l}(y, \mu) \sim& \cosmu(\mu) \sum_{w\in W} C_1(\mu^w) p_{\rho+\mu^w}(y) \int_{Y^+} W^*(t,-\mu,\psi_{1,1}) W(t,\mu^w,\psi_{1,1}) \, t_1^2 t_2 \, dt \\
	=& \frac{\cosmu(\mu)}{2\sqrt{\pi}} \sum_{w\in W} \frac{C_1(\mu^w)}{\Lambda(\mu^w)} p_{\rho+\mu^w}(y) \prod_{j,k} \Gamma\paren{\frac{1+\mu_j-\mu_k}{2}},
\end{align*}
by Stade's formula \eqref{eq:StadesFormula}.
Note our first-term asymptotic is independent of the sign of $y$; this is because all of the sign information in $y$ is carried by the term $\psi_y(x^*)$ which drops out in the limit.

If we can justify the above analysis, then we must have
\begin{align}
\label{eq:KwlValue}
	K_{w_l}(y, \mu) =& -\frac{\pi^3}{32} \sum_{w\in W} \frac{1}{\sinmu(\mu^w)} J_{w_l}(y,\mu^w),
\end{align}
after simplifiying the gamma factors.
Note that this does not imply an equality for $K_{w_l}^\text{as}$.

\subsection{Justification of the long element term}
\label{sect:LongEleJustify}
As in \cite{Me01}, consider
\begin{align}
\label{eq:LEWeightFuncDelta}
	H_{w_l}^*(f; y, s) =& \int_{G/K} \int_{\Re(\mu)=0} W^*(z,-\mu,\psi_{1,1}) W^*(y w_l z,-\mu,\psi_{1,1}) \\
	& \qquad f(\mu) \, \specmu(\mu) d\mu \, p_{(s,s,-2s)}(z) \, dz, \nonumber
\end{align}
with $s \ge 1$.
This extends the previous definition slightly with
\[ H_{w_l}(f; y) = \frac{1}{96 \pi^6 \abs{y_1 y_2}} H_{w_l}^*(f; y, 1). \]
We justify the previous section in the following steps:
\begin{enumerate}
\item As in section \ref{sect:HardTerms}, we may construct $K_{w_l}^\text{as}(y,\mu,s)$ so that
\begin{align}
\label{eq:HwusingKwas}
	H_{w_l}^*(f; y, s) = \int_{\Re(\mu)=(\eta,0,-\eta)} f(\mu) K_{w_l}^\text{as}(y, \mu, s) d\mu,
\end{align}
for some small $\eta > 0$, where $K_{w_l}^\text{as}(y, \mu, s)$ is defined by absolutely convergent integrals on $\Re(\mu_1) > \Re(\mu_2) > \Re(\mu_3)$.

\item \label{eq:JustifyAnalCont} Starting with $\Re(s) \ge 1000$, we will show $K_{w_l}^\text{as}(y, \mu, s)$ extends holomorphically in $\mu$ to
\[ \Re(\mu_1-\mu_2),\Re(\mu_2-\mu_3) > -\eta, \]
and is twice differentiable in each $y_i$ there.

\item \label{item:JustifyNonsense} Again, we extend both $H_{w_l}^*$ and $K_{w_l}^\text{as}$ to functions on $G$, and we will show that
\[ K_{w_l}(g, \mu, s) := \frac{1}{6} \sum_{w\in W} K_{w_l}^\text{as}(g,\mu^w,s) \]
satisfies the expected differential equations.
This implies, by our explicit solution of these differential equations in section \ref{sect:LongEleDiffEqs},
\[ K_{w_l}(y, \mu, s) = \sum_{w\in W} C_2(\mu^w, s) J_{w_l}(y,\mu^w), \]
for some function $C_2(\mu,s)$ independent of $y$.

\item \label{eq:JustifyAsymp} We will justify the first-term asymptotic of $K_{w_l}(y,\mu,s)$ at large $s$, which gives
\begin{align*}
	C_2(\mu,s) =& \frac{\Gamma\paren{\frac{s}{2}}^3}{32\pi^{3/2}} \frac{\cosmu(\mu)}{\sinmu(\mu)} \prod_{j<k} \Gamma\paren{\frac{s+\mu_j-\mu_k}{2}} \Gamma\paren{\frac{s+\mu_k-\mu_j}{2}},
\end{align*}
and, as in section \ref{sect:AsympAnal}, we have the particular value
\begin{align*}
	C_2(\mu,1) =& \frac{\pi^3}{32 \sinmu(\mu)}.
\end{align*}

\item At last, we have analytic continuation to $\Re(s) > \eta$ of both the original integral \eqref{eq:LEWeightFuncDelta} (which can be seen through \eqref{eq:HwusingKwas}), and the final form
\[ H_{w_l}^*(f; y, s) = \int_{\Re(\mu)=0} f(\mu) K_{w_l}(y, \mu, s) d\mu, \]
so we may take $s=1$, and this completes the construction.
\end{enumerate}
We will give the steps slightly out of order, showing \ref{item:JustifyNonsense}, then \ref{eq:JustifyAnalCont}, and finally \ref{eq:JustifyAsymp}.
It should be mentioned at this point that we will only require some polynomial bound in $\mu$ on the analytic continuation of $K_{w_l}^\text{as}$ for step \ref{item:JustifyNonsense} and for step \ref{eq:JustifyAsymp}, we may even assume that $\mu$ is fixed.

\subsubsection{The differential equations}
We first show step \ref{item:JustifyNonsense} above, assuming the continuation and differentiability of $K_{w_l}^\text{as}$.
For $f(\mu)$ holomorphic on a neighborhood of $\Re(\mu)=0$, symmetric in $\mu$, and of rapid decay in each $\abs{\Im(\mu_i)}$, consider the integral
\begin{align}
\label{eq:HiDef}
	H_i :=& \int_{\Re(\mu)=0} f(\mu) \, \wtilde{\Delta_i} K_{w_l}^\text{as}(y, \mu, s) \, \specmu(\mu) d\mu,
\end{align}
for $i=0,1$.
By definition of $\wtilde{\Delta_i}$ and $K_{w_l}^\text{as}(g, \mu, s)$, and by absolute and uniform convergence of the $\mu$ integral, we have
\begin{align*}
	H_i =& \left. \Delta_i \int_{\Re(\mu)=0} f(\mu) K_{w_l}^\text{as}(g, \mu, s) \, \specmu(\mu) d\mu \right|_{u=v=I}.
\end{align*}
Of course, by construction of $K_{w_l}^\text{as}$, this is
\begin{align*}
	H_i =& \left. \Delta_i H_{w_l}^*(f; g, s) \right|_{u=v=I},
\end{align*}
but before reconstructing the Whittaker function from the integral representation \eqref{eq:WhittakerMellinShifted}, we could have brought the differential operator inside the $z$ integral
\begin{align*}
	H_i =& \int_{G/K} \Delta_i \int_{\Re(\mu)=0} W^*(z,-\mu,\psi_{1,1}) W^*(g w_l z,-\mu,\psi_{1,1}) \\
	& \qquad \left. f(\mu) \, \specmu(\mu) d\mu \, p_{(s,s,-2s)}(z) \, dz \right|_{u=v=I}, \nonumber
\end{align*}
by absolute and uniform convergence of the $z$ integral (via contour shifting in $\mu$).
Now we may finally bring $\Delta_i$ inside the $\mu$ integral, apply it to the Whittaker function, and reapply the construction of $K_{w_l}^\text{as}$,
\begin{align}
\label{eq:HiResult}
	H_i =& \int_{\Re(\mu)=0} f(\mu) \, \lambda_i(\mu) \, K_{w_l}^\text{as}(y, \mu, s) \, \specmu(\mu) d\mu.
\end{align}

Combining \eqref{eq:HiDef} and \eqref{eq:HiResult}, we have that $K_{w_l}$ satisfies the appropriate differential equations in the sense of distributions:
\begin{align*}
	\int_{\Re(\mu)=0} f(\mu) \, \wtilde{\Delta_i} K_{w_l}^\text{as}(y, \mu, s) \, \specmu(\mu) d\mu =& \int_{\Re(\mu)=0} f(\mu) \, \lambda_i(\mu) \, K_{w_l}^\text{as}(y, \mu, s) \, \specmu(\mu) d\mu.
\end{align*}
To remove the $\mu$ integral, we may take the limit $\delta \to 0$ with
\[ f(\mu) = f(\mu,\mu',\delta) = \sum_{w\in W} \frac{1}{\delta^2} \exp\paren{\frac{1}{\delta^2}\sum_{i=1}^3(\mu_i^w-\mu'_i)^2}, \]
giving
\[ \wtilde{\Delta_i} \sum_{w\in W} K_{w_l}^\text{as}(y, {\mu'}^w, s) = \lambda_i(\mu') \sum_{w\in W} K_{w_l}^\text{as}(y, {\mu'}^w, s), \]
and the previous analysis of these differential equations gives
\[ K_{w_l}(y, \mu, s) = \sum_{w\in W} C_2(\mu^w, s) J_{w_l}(y,\mu^w). \]
We should point out here that this is the step which would not have worked away from $\Re(\mu)=0$.

\subsubsection{The analytic continuation}
\label{sect:AnalCont}
We now show the analytic continuation of step \ref{eq:JustifyAnalCont} above.
We do this by writing $T_{w_l}(y,\mu,\mu',s)$, appropriately generalized from \eqref{eq:TwDef}, in a form that converges in a neighborhood of $\Re(\mu)=0$.
The form of our answer is sufficiently complicated that we will not write it out, but rather describe its construction.

First, we isolate the $x'$ integral,
\begin{align}
\label{eq:TwADef}
	T_{w_l}(y, \mu, \mu') =& \frac{\cosmu(\mu)}{4\pi^2} p_{\rho+\mu'}(y \, t^{w_l}) \int_{Y^+} W^*(t,-\mu, \psi_{1,1}) X'(y, t, \mu, \mu') \, t_1^{2+2s} t_2^{2+s} \,dt. \\
\label{eq:XprimeDef}
	X'(y, t, \mu, \mu') =& \int_{\wbar{U}_{w_l}(\R)} p_{\rho+\mu'}(y^*) \, \psi_{y \, t^{w_l}}(x^*) \psi_t(x) dx.
\end{align}
Here we have substituted $x't \mapsto tx$, dropping the prime notation for convenience; in particular this means that
\[ x^* y^* \equiv w_l x \pmod{\R^+ K}. \]

The exponential decay of the Whittaker function in large $t_i$, and the high power of $t_i$ coming from the term $(t_1^2 t_2)^s$ with $s$ large means that we may treat each $t_i \approx 1$.
The asymptotics for large $t_i$ have some cost in terms of $\mu$, but we will only require polynomial dependence on $\mu$ in steps \ref{item:JustifyNonsense} and \ref{eq:JustifyAsymp} of the justification.
In other words, we will analyze just the $X'$ integral, ignoring the positive and negative powers of $t_i$ that result from integration by parts, as this will not affect the convergence of the whole $T_{w_l}$ integral, and we will not explicitly track the dependence on $\mu$.

Our computations require the explicit form of $X'$.
First, the components of $x^*$ and $y^*$ are
\begin{align}
\label{eq:xstar}
	x^*_1 =& -\frac{x_2+x_1x_3}{1+x_2^2+x_3^2}, & x^*_2 =& -\frac{x_1+x_2(x_1x_2-x_3)}{1+x_1^2+(x_1x_2-x_3)^2}, \\
\label{eq:ystar}
	y^*_1 =& \frac{\sqrt{1+x_1^2+(x_1x_2-x_3)^2}}{1+x_2^2+x_3^2}, & y^*_2 =& \frac{\sqrt{1+x_2^2+x_3^2}}{1+x_1^2+(x_1x_2-x_3)^2}.
\end{align}
Then separating the $x_3$ integral, we have
\begin{align*}
	X'(y, t, \mu, \mu') =& \int_{\R^2} \e{t_1 x_1+t_2 x_2} X'_3\paren{x_1, x_2,\frac{y_1}{t_2},\frac{y_2}{t_1},s_1,s_2,0,0} dx, \\
	X'_3(x_1, x_2,\alpha_1,\alpha_2,s_1,s_2,n_1,n_2) =& \int_\R (x_1x_2-x_3)^{n_1} x_3^{n_2} (1+x_1^2+(x_3-x_1x_2)^2)^{s_1} (1+x_2^2+x_3^2)^{s_2} \\
	& \qquad \e{-\alpha_1\frac{x_2+x_1 x_3}{1+x_2^2+x_3^2}-\alpha_2\frac{x_1+x_2(x_1x_2-x_3)}{1+x_1^2+(x_1 x_2-x_3)^2}} dx_3,
\end{align*}
where $s_1=-\frac{1+\mu'_1-\mu'_2}{2}$ and $s_2=-\frac{1+\mu'_2-\mu'_3}{2}$, and for convenience, we write
\[ X'_3(s_1,s_2,n_1,n_2) = X'_3(x_1, x_2,\alpha_1,\alpha_2,s_1,s_2,n_1,n_2). \]
To achieve the holomorphy requirements, we need to rearrange the $X'$ integral such that it converges absolutely on $\Re(\mu_1'-\mu_2'),\Re(\mu_2'-\mu_3')>-\eta$.

Now the integral over the region $\Max{\abs{x_1},\abs{x_2}} < 1$ already converges on $\Re(\mu)=0$, so we need only consider the regions $\abs{x_2} > \Max{1,\abs{x_1}}$ and $\abs{x_1} > \Max{1,\abs{x_2}}$.
We address the first range, the second is identical by symmetry (after sending $x_3 \mapsto x_3+x_1x_2$).
If $M_1(x)$ is a smooth function which is 1 on $\abs{x} < \frac{1}{2}$ and 0 on $\abs{x} > 2$, we may select $\abs{x_2} > \abs{x_1}$, up to a constant, by inserting the function
\[ M_2(x) = \frac{1}{2}\paren{1+M_1(x)-M_1(1/x)}, \]
which satisfies $M_2(x)=1-M_2(1/x)$; thus both $M_2\paren{\frac{1+x_1^2}{1+x_2^2}}$ and the complementary function $1-M_2\paren{\frac{1+x_1^2}{1+x_2^2}}=M_2\paren{\frac{1+x_2^2}{1+x_1^2}}$ save powers of the larger variable under differentiation on the appropriate range.
Then we may isolate $\abs{x_2} > 1$ by inserting a function $M_1^c(\sqrt{\abs{x_2}})$, where $M_1^c=1-M_1$, whose derivatives save powers of $x_2$.
Note that we do not need the complementary function $M_1(\sqrt{\abs{x_2}})$ to be differentiable at 0 since we need do no integrating by parts on the range $\abs{x_1}<\abs{x_2}<1$.

Now let $\eta$ be a fixed, sufficiently small, positive real number.
For the desired convergence, we need only acquire a factor $\abs{x_2}^{-10\eta}$ through integration by parts.
The strategy is to integrate by parts in $x_2$, which solves the problem outside a certain problematic region, on which we integrate by parts in $x_3$.
The final function constructed will be polynomially bounded in $t, t^{-1}$ and $\Im(\mu')$, as required, but we will not track this explicitly.

The method of \cite[sect. 4.3]{Me01} works again here:
Substituting $x_3 \mapsto x_2 x_3$ gives
\begin{align*}
	X'_3(s_1,s_2,n_1,n_2) =& x_2^{n_1+n_2+1+2(s_1+s_2)} \int_\R (x_1-x_3)^{n_1} x_3^{n_2} (1+x_2^{-2}+x_3^2)^{s_1} \paren{\frac{1+x_1^2}{x_2^2}+(x_3-x_1)^2}^{s_1} \\
	& \qquad  \e{-\alpha_1 \frac{1}{x_2}\frac{1+x_1 x_3}{1+x_2^{-2}+x_3^2}-\alpha_2\frac{x_1 x_2^{-2}+(x_1-x_3)}{\frac{1+x_1^2}{x_2^2}+(x_1-x_3)^2}} dx_3,
\end{align*}
so
\begin{align}
\label{eq:x2derv}
	& \pderv{}{x_2} X'_3(s_1,s_2,n_1,n_2) = \\
	& \qquad \frac{n_1+n_2+1+2(s_1+s_2)}{x_2} X'_3(s_1,s_2,n_1,n_2) -2\frac{s_2}{x_2}X'_3(s_1,s_2-1,n_1,n_2) \nonumber\\
	& \qquad -2s_1\frac{1+x_1^2}{x_2}X'_3(s_1-1,s_2,n_1,n_2) +\alpha_1 X'_3(s_1,s_2-1,n_1,n_2) \nonumber\\
	& \qquad - 2\alpha_1 x_2 X'_3(s_1,s_2-2,n_1,n_2) +\alpha_1 \frac{x_1}{x_2} X'_3(s_1,s_2-1,n_1,n_2+1) \nonumber\\
	& \qquad -2\alpha_1 x_1 X'_3(s_1,s_2-2,n_1,n_2+1) +2\alpha_2 \frac{x_1}{x_2}X'_3(s_1-1,s_2,n_1,n_2) \nonumber\\
	& \qquad +2\alpha_2 \frac{x_1(1+x_1^2)}{x_2} X'_3(s_1-2,s_2,n_1,n_2) +2\alpha_2 (1+x_1^2) X'_3(s_1-2,s_2,n_1+1,n_2).\nonumber
\end{align}

Now in integrating by parts in $x_2$, only the last term of \eqref{eq:x2derv} fails to save powers of $x_2$, and only in the range
\[ \abs{x_1} < \abs{x_2}^\eta, \qquad \abs{x_1x_2-x_3} < \abs{x_2}^{2\eta}. \]
Selecting this range by a smooth partition of unity, we may expand most of the terms of the complicated exponential in the first few terms of their power series; the remainder terms will satisfy the holomorphy requirements. This gives
\begin{align*}
	X' =& \int_{\R^3} M_3(x) (1+x_1^2+(x_3-x_1x_2)^2)^{s_1} (1+x_2^2+x_3^2)^{s_2} \\
	& \qquad \e{t_1 x_1+t_2 x_2-\alpha_2 \frac{x_1+x_2(x_1x_2-x_3)}{1+x_1^2+(x_1 x_2-x_3)^2}} dx+\ldots,
\end{align*}
where the $\ldots$ indicates terms which are either lower-order or already meet the holomorphy requirements, and
\[ M_3(x) = M_2\paren{\tfrac{1+x_1^2}{1+x_2^2}} M_1^c(\sqrt{\abs{x_2}}) \, M_1\paren{\tfrac{\paren{1+x_1^2}^{1/\eta}}{x_2^2}} M_1\paren{\tfrac{\paren{1+(x_1 x_2-x_3)^2}^{1/\eta}}{x_2^4}}. \]
Now send $x_3 \mapsto x_3\sqrt{1+x_1^2}+x_1x_2$,
\begin{align*}
	X' =& \int_{\R^3} M_2(x) (1+x_1^2)^{s_1+\frac{1}{2}} (1+x_3^2)^{s_1} \paren{1+x_2^2+\paren{x_3\sqrt{1+x_1^2}+x_1x_2}^2}^{s_2} \\
	& \qquad \e{t_1 x_1+t_2 x_2+\alpha_2 \frac{x_1}{(1+x_1^2)(1+x_3^2)}+\alpha_2 \frac{x_2}{\sqrt{1+x_1^2}}\frac{x_3}{1+x_3^2}} dx+\ldots,
\end{align*}
where
\[ M_4(x) = M_3\paren{x_1, x_2, x_3\sqrt{1+x_1^2}+x_1x_2}. \]

We need to exclude the three problem points $x_3=0,\pm 1$:
First, for $\abs{x_3} < \abs{x_2}^{-6\eta}$, we have
\[ \frac{x_3}{1+x_3^2} = x_3-x_3^3+\ldots, \]
so by truncating at roughly $\frac{1}{2\eta}$ terms inside the exponential, the remainder term meets the holomorphy and differentiability requirements.
Then integrating by parts with the lead term,
\[ \e{\alpha_2 \frac{x_2}{(1+x_1^2)}x_3}, \]
each iteration saves at least $\abs{x_2}^{10\eta}$.
The same trick applies for $x_3=1+u$ with $\abs{u}<\abs{x_2}^{-6\eta}$, where
\[ \frac{x_3}{1+x_3^2} =\frac{1}{2}-\frac{1}{4} u^2+\ldots, \]
and we integrate by parts with
\[ u \, \e{\alpha_2 \frac{x_2}{4(1+x_1^2)}u^2}, \]
(the next term in the $u$ expansion is $\BigO{u^3}$, so we gain $\frac{1}{u} \derv{}{u}u^3=\BigO{\abs{x_2}^{-6\eta}}$), and we repeat the trick again for $x_3$ near $-1$.
Otherwise, we may send $\frac{x_3}{1+x_3^2} \mapsto u$, giving
\[ x_3 = \frac{1\pm\sqrt{1-u^2}}{u}, \qquad \derv[n]{x_3}{u} \ll \paren{\abs{x_2}^{20\eta}}^n, \]
and integrate by parts in $u$, with each iteration saving at least $\abs{x_2}^{1-30\eta}$.
This completes the construction.

\subsubsection{The differentiability}
We now show the differentiability of step \ref{eq:JustifyAnalCont} above.
It is trivial to see from \eqref{eq:xstar} that
\[ \abs{x_1^*} < 1+\abs{x_1}, \qquad \abs{x_2^*} < 1 + \abs{x_2}, \]
so we will have differentiability of the continued function at $\Re(\mu)=0$, provided we can save somewhat more, say $\abs{x_2}^5$ on the range $\abs{x_2} > \Max{1,\abs{x_1}}$, in the integration by parts, but this is no more difficult than saving the $\abs{x_2}^{10\eta}$ we already did.
Thus the differentiability follows simply by more repetitions of the integration by parts.

\subsubsection{The asymptotics}
At last, we show step \ref{eq:JustifyAsymp} above.
This can be accomplished by Mellin expanding the complicated exponential term as was done in \cite[sect. 5]{Me01} and shifting integral, but this is arguably even less pleasant than the integration by parts argument of section \ref{sect:AnalCont}.
Instead, we attack the problem by writing
\[ \psi_{y \, t^{w_l}}(x^*) = \paren{\psi_{y \, t^{w_l}}(x^*)-1}+1 \]
in \eqref{eq:XprimeDef}.
Then, as in the asymptotic analysis of section \ref{sect:LongEleTerm}, the $X'$ integral without the complicated exponential becomes the Jacquet integral of the Whittaker function,
\begin{align*}
	X'(y, t, \mu, \mu') =& (t_1 t_2)^{-2} p_{-\rho-\mu'}(t^{w_l}) W(t,\mu',\psi_{1,1})+\wtilde{X'}(y,t,\mu,\mu') \\
	\wtilde{X'}(y,t,\mu,\mu') =& \int_{\wbar{U}_w(\R)} p_{\rho+\mu'}(y^*) \, (\psi_{y \, t^{w_l}}(x^*)-1) \psi_t(x) dx.
\end{align*}
The difference from the argument in the earlier section is that we are now in the region of absolute convergence of the Jacquet integral.

Now the analytic continuation, differentiability, and asymptotics of the Whittaker function are known, and we apply the integration by parts argument of section \ref{sect:AnalCont} to $\wtilde{X'}$.
If we can show that the construction, when applied to $\wtilde{X'}$, produces something which is small in $y$ as $y \to 0$, with polynomial dependence on $\mu$ and $t_i^{\pm 1}$, then the remainder of the asymptotic analysis of section \ref{sect:AsympAnal} becomes rigorous.
As before, the rapid decay in the $T_w$ integral means we may assume that
\begin{align}
\label{eq:Xprimetrange}
	t_1,t_2,t_1^{-1},t_2^{-1} < \Min{\abs{y_1}^{-\eta},\abs{y_2}^{-\eta}}.
\end{align}

First, we remove the trivial case:
On the region $\Max{\abs{x_1}, \abs{x_2}} < 1$, we have $\abs{x_1^*},\abs{x_2^*} < 2$, so 
\[ \abs{\psi_{y \, t^{w_l}}(x^*)-1} \ll \abs{y_1}^{1-\eta} + \abs{y_2}^{1-\eta}. \]

Next, we examine the region $\abs{x_2} > \Max{\abs{x_1},1}$:
Differentiating the term $\psi_{y \, t^{w_l}}(x^*)-1$ naturally produces positive powers of $y_1$ or $y_2$, so we may ignore those terms, except:
The integration by parts in $x_3$ produces negative powers of $\alpha_2 = \frac{y_2}{t_1}$, but we only do this for the final term of \eqref{eq:x2derv}, which itself produces positive powers of $\alpha_2$.
Hence we may always arrange for the net power of $y_2$ to be positive by simply integrating by parts in $x_2$ more times than in $x_3$.

At last, we must consider those terms in the integration by parts construction which do not result from differentiating $\psi_{y \, t^{w_l}}(x^*)-1$; these are essentially derived the first three terms on the right-hand side of \eqref{eq:x2derv}.
A typical term is of the form
\begin{align*}
	& P(s,n) \int_{U(\R)} M_2\paren{\frac{1+x_1^2}{1+x_2^2}} M_1^c(\sqrt{\abs{x_2}}) (x_1 x_2-x_3)^{n_1} x_3^{n_2} x_2^{-n_3} (1+x_1^2)^{n_4} \\
	& \qquad (1+x_1^2+(x_3-x_1 x_2)^2)^{s_1} (1+x_2^2+x_3^2)^{s_2} \paren{\psi_{y \, t^{w_l}}(x^*)-1} dx,
\end{align*}
where $P(s,n)$ is some polynomial in $s$ and $n$, and
\[ n_1+n_2+2\Re(s_1+s_2) < -1, \qquad n_1+n_2-n_3+2n_4+2\Re(s_1+s_2) \le -6, \] 
which is a technical way of saying the $x_3$ integral converges and the $x_2$ integral converges quickly.
The argument for the region $\Max{\abs{x_1}, \abs{x_2}} < 1$ is now essentially sufficient, due to the overconvergence of the $x_2$ integral:
The portion of the integral with $\abs{x_2} > \Max{\abs{y_1}^{-1/4}, \abs{y_2}^{-1/4}}$ is trivially bounded by $\abs{y_1}^{1/2}+\abs{y_2}^{1/2}$, and the remainder is bounded by $\abs{y_1}^{1/2-\eta}+\abs{y_2}^{1/2-\eta}$ using
\[ \abs{\psi_{y \, t^{w_l}}(x^*)-1} \ll \abs{y_1}^{1-\eta}\abs{x_1^*}+\abs{y_2}^{1-\eta}\abs{x_2^*} \ll \abs{y_1}^{3/4-\eta}+\abs{y_2}^{3/4-\eta}. \]

The computation of $C_2(\mu,s)$ now proceeds as in section \ref{sect:AsympAnal}, using the explicit asymptotic
\begin{align*}
	X'(y, t, \mu, \mu') =& (t_1 t_2)^{-2} p_{-\rho-\mu'}(t^{w_l}) W(t,\mu',\psi_{1,1})+O_{t,\mu}\paren{\abs{y_1}^{1/2-\eta}+\abs{y_2}^{1/2-\eta}},
\end{align*}
on the region \eqref{eq:Xprimetrange}, where the implied constant depends polynomially on $\mu$, $t$, and $t^{-1}$.
The parameter $s$ enters the computations through Stade's formula \eqref{eq:StadesFormula}.

We note that it is possible to give an expression for $X'$ which simultaneously shows the analytic continuation, differentiability and asymptotics, but this would be quite lengthy.

\subsection{\texorpdfstring{The $\boldsymbol{w_4}$ term}{The w4 term}}
\label{sect:w4Term}
As in section \ref{sect:LongEleTerm}, we consider $K_{w_4}(g,\mu)$ to be defined by the non-convergent integral \eqref{eq:JwFakeDef},
\begin{align}
\label{eq:Kw4FakeDef}
	K_{w_4}(g,\mu) =& \int_{Y^+} \int_{\wbar{U}_{w_4}(\R)} W\paren{x t,-\mu, \psi_{1,1}} W\paren{g w_4 x t, \mu, \psi_{1,1}}  dx \, t_1^2 t_2 \, dt,
\end{align}
and the $x$ integral here is over the space
\[ \wbar{U}_{w_4}(\R)=\set{\Matrix{1&x_2&x_3\\&1&0\\&&1}\setdiv x_2,x_3\in\R}. \]

\subsubsection{The differential equations}
As before, we may substitute
\[ w_4^{-1} \Matrix{1\\0&1\\v_3&v_2&1} w_4 \, x = \Matrix{1&v_3&v_2\\&1&0\\&&1}x \mapsto x, \]
to see that
\[ K_{w_4}\paren{ruy\SmallMatrix{1\\v_1&1\\v_3&v_2&1},\mu} = \e{u_1+u_2-v_3} K_{w_4}\paren{y\SmallMatrix{1\\v_1&1\\0&0&1},\mu}, \]
so the partial derivatives in $u$ and $v$ are
\[ \partial_{u_1} K_{w_4}(g,\mu) = \partial_{u_2} K_{w_4}(g,\mu) = -\partial_{v_3} K_{w_4}(g,\mu) = 2\pi i K_{w_4}(g,\mu), \]
and
\[ \partial_{u_3} K_{w_4}(g,\mu) = \partial_{v_2} K_{w_4}(g,\mu) = 0. \]
We don't have a good way of determining the dependence on $v_1$, but the particular combination of operators given by \eqref{eq:Delta3op} will avoid this, because all of the derivatives in $v_1$ drop out.
We have $\Delta_3 \, K_{w_4}(g,\mu) = 0$, and this reduces to
\begin{align}
\label{eq:Kw4Delta3}
	\what{\Delta_3} \, K_{w_4}(y,\mu) = 0,
\end{align}
where
\begin{align}
\label{eq:Kw4Delta3Def}
	\what{\Delta_3} = \lambda_1(\mu)+\lambda_2(\mu) + 8\pi^3 i y_1 y_2-\lambda_1(\mu) y_1 \partial_{y_1}-y_1^3 \partial_{y_1}^3.
\end{align}
For a fixed $y_2$, we may solve the differential equation in $y_1$ with a power series as before, and the solutions are given by
\[ J_{w_4}(y,\mu) = \abs{8\pi^3 y_1}^{1-\mu_3} \sum_{n=0}^\infty \frac{(8\pi^3 i y_1 y_2)^n}{n! \, \Gamma\paren{n+1+\mu_1-\mu_3} \, \Gamma\paren{n+1+\mu_2-\mu_3}}, \]
along with $J_{w_4}(y,\mu^{w_4})$ and $J_{w_4}(y,\mu^{w_5})$.
It follows that
\begin{align}
\label{eq:Kw4Value}
	K_{w_4}(y,\mu) = \sum_{w\in W_3} C_3(\mu^w, y_2) J_{w_4}(y,\mu^w),
\end{align}
for some function $C_3(\mu,y_2)$ which is independent of $y_1$.

\subsubsection{The asymptotic analysis}
\label{sect:MultConst}
To determine the value of $C_3(\mu,-1)$, we take the first-term asymptotic in $y_1$ in \eqref{eq:Kw4FakeDef} at $u=v=I$.
Since we are only taking the limit in one of the $y$ variables, the resulting expression is much more complicated than for the long element, but this can still be evaluated using the Barnes integral lemmas.

We apply \eqref{eq:Whitty1Asymp}, the first-term asymptotic in $y_1$ of the Whittaker function, so that
\begin{align*}
	K_{w_4}(y,\mu) \sim& \frac{\cosmu(\mu)}{8\pi^2} \sum_{w\in W_3} \abs{\pi y_1}^{1-\mu_3^w} \abs{\pi y_2}^{1-\frac{1}{2}\mu_3^w} \Gamma\paren{\tfrac{\mu_3^w-\mu_1^w}{2}} \Gamma\paren{\tfrac{\mu_3^w-\mu_2^w}{2}} \int_0^\infty \int_0^\infty W^*\paren{t,-\mu, \psi_{1,1}} \\
	& \int_{\wbar{U}_{w_4}(\R)} \e{t_2 x_2+y_2 t_1 x_2^*} (y_1^* y_2^*) ((y_1^*)^2 y_2^*)^{-\frac{1}{2}\mu_3^w} K_{\frac{\mu_1^w-\mu_2^w}{2}}\paren{2\pi \abs{y_2} t_1 y_2^*} dx \, (t_1 t_2^2)^{\frac{1}{2}\mu_3^w} \frac{dt_1 dt_2}{t_2},
\end{align*}
where $x^* y^* \equiv w_4 x \pmod{\R^+ SO(3,\R)}$.
Explicitly,
\[ x_2^* = -\frac{x_2 x_3}{1+x_2^2}, \quad y_1^* = \frac{\sqrt{1+x_2^2}}{1+x_2^2+x_3^2}, \quad y_2^* = \frac{\sqrt{1+x_2^2+x_3^2}}{1+x_2^2}. \]

We substitute $x_3 \mapsto x_3 \sqrt{1+x_2^2}$ and Mellin expand the Whittaker and $K$-Bessel functions using \eqref{eq:WhittakerMellin} and \eqref{eq:KBesselMellin},
\begin{align*}
	K_{w_4}(y,\mu) \sim& \frac{\cosmu(\mu)}{32\pi^4} \sum_{w\in W_3} \abs{\pi y_1}^{1-\mu_3^w} \abs{\pi y_2}^{1-\frac{1}{2}\mu_3^w} \Gamma\paren{\tfrac{\mu_3^w-\mu_1^w}{2}} \Gamma\paren{\tfrac{\mu_3^w-\mu_2^w}{2}} \int_0^\infty \int_0^\infty \int_0^\infty \int_0^\infty \\
	& \qquad \int_{\Re(s)=(\frac{3}{2},\frac{1}{4},\frac{1}{4})} G((s_1,s_2),-\mu) \Gamma\paren{\tfrac{2s_3+\mu_1^w-\mu_2^w}{4}} \Gamma\paren{\tfrac{2s_3+\mu_2^w-\mu_1^w}{4}} \\
	& \qquad (\pi t_1)^{1-s_1-s_3} (\pi t_2)^{1-s_2} \abs{y_2}^{-s_3} (1+x_2^2)^{\frac{-2+2s_3+3\mu_3^w}{4}} (1+x_3^2)^{\frac{-2-2s_3+3\mu_3^w}{4}} \\
	& \qquad \cos\paren{2\pi t_2 x_2}\cos\paren{2\pi y_2 t_1 \tfrac{x_2 x_3}{\sqrt{1+x_2^2}}} \frac{ds}{(2\pi i)^3} dx_2 \, dx_3 \, (t_1 t_2^2)^{\frac{1}{2}\mu_3^w} \frac{dt_1 \,dt_2}{t_2}.
\end{align*}

Now we may evaluate the $t$ and then $x$ integrals using \eqref{eq:CosMellin} and \eqref{eq:BetaDef}
\begin{align*}
	K_{w_4}(y,\mu) \sim& \frac{\cosmu(\mu)}{128\pi^4} \sum_{w\in W_3} 2^{-3-\frac{3}{2}\mu_3^w} \pi^{1-3\mu_3^w} \abs{y_1}^{1-\mu_3^w} \abs{y_2}^{-1-\mu_3^w} \Gamma\paren{\tfrac{\mu_3^w-\mu_1^w}{2}} \Gamma\paren{\tfrac{\mu_3^w-\mu_2^w}{2}} \\
	& \int_{\Re(s)=(\frac{3}{2},\frac{1}{4},\frac{1}{4})} G((s_1,s_2),-\mu) \Gamma\paren{\tfrac{2s_3+\mu_1^w-\mu_2^w}{4}} \Gamma\paren{\tfrac{2s_3+\mu_2^w-\mu_1^w}{4}} 2^{s_1+s_2+s_3} \abs{y_2}^{s_1} \\
	& \Gamma\paren{1-s_2+\mu_3^w} \cos\frac{\pi}{2}\paren{1-s_2+\mu_3^w} \Gamma\paren{2-s_1-s_3+\tfrac{1}{2}\mu_3^w} \cos\frac{\pi}{2}\paren{2-s_1-s_3+\tfrac{1}{2}\mu_3^w} \\
	& B\paren{\tfrac{-4+2s_1+2s_2+2s_3-3\mu_3^w}{4},\tfrac{2-2s_2-2s_3-\mu_3^w}{4}} B\paren{\tfrac{-2+2s_1+2s_3-\mu_3^w}{4},\tfrac{2-s_1-\mu_3^w}{2}} \frac{ds}{(2\pi i)^3}.
\end{align*}
The interchange of integrals here is much easier to justify than those of previous sections; alternately, instead of applying the Mellin expansion of the Whittaker function, one may consider the Mellin expansion of the cosines as in \cite[eq. (31)]{Me01}.

We apply the reflection and duplication formulas \eqref{eq:ReflectionDuplication} to the cosines and gamma functions, respectively, on the third line of the previous display.
Then the $s_2$ integral of just the gamma functions which depend on $s_2$ becomes
\begin{align*}
	&\int_{\Re(s_2)=\frac{1}{4}} \frac{\Gamma\paren{\frac{s_2-\mu_1^w}{2}} \Gamma\paren{\frac{s_2-\mu_2^w}{2}} \Gamma\paren{\frac{1-s_2+\mu_3^w}{2}} \Gamma\paren{\tfrac{-4+2s_1+2s_2+2s_3-3\mu_3^w}{4}} \Gamma\paren{\tfrac{2-2s_2-2s_3-\mu_3^w}{4}}}{\Gamma\paren{\frac{s_1+s_2}{2}}} \frac{ds_2}{2\pi i} \\
	&= 2 
\frac{\Gamma\paren{\frac{1+\mu_3^w-\mu_1^w}{2}} \Gamma\paren{\frac{1+\mu_3^w-\mu_2^w}{2}} \Gamma\paren{\frac{-2+2s_1+2s_3-\mu_3^w}{4}} \Gamma\paren{\tfrac{2-2s_3+\mu_2^w-\mu_1^w}{4}} \Gamma\paren{\tfrac{2-2s_3+\mu_1^w-\mu_2^w}{4}} \Gamma\paren{\tfrac{-1+s_1-2\mu_3^w}{2}}}{ \Gamma\paren{\frac{s_1+\mu_1^w}{2}} \Gamma\paren{\frac{s_1+\mu_2^w}{2}} \Gamma\paren{\frac{4-2s_3+3\mu_3^w}{4}}},
\end{align*}
by the second Barnes integral lemma, Theorem \ref{thm:BarnesSecond}.

Now at $\abs{y_2}=1$, the $s_1$ integral of just those gamma functions which depend on $s_1$ becomes
\begin{align*}
	& \int_{\Re(s_1)=\frac{3}{2}} \Gamma\paren{\tfrac{s_1+\mu_3^w}{2}} \Gamma\paren{\tfrac{2-s_1-\mu_3^w}{2}} \Gamma\paren{\tfrac{4-2s_1-2s_3+\mu_3^w}{4}} \Gamma\paren{\tfrac{-2+2s_1+2s_3-\mu_3^w}{4}} \frac{ds_1}{2\pi i} \\
	&= 4 \, \Gamma\paren{\tfrac{2+2s_3-3\mu_3^w}{4}} \Gamma\paren{\tfrac{4-2s_3+3\mu_3^w}{4}},
\end{align*}
by the first Barnes integral lemma, Theorem \ref{thm:BarnesFirst}.

At last, the $s_3$ integral of just those gamma functions depending on $s_3$ may be evaluated by the first Barnes integral lemma,
\begin{align*}
	& \int_{\Re(s_3)=\frac{1}{4}} \Gamma\paren{\tfrac{2s_3+\mu_1^w-\mu_2^w}{4}} \Gamma\paren{\tfrac{2s_3+\mu_2^w-\mu_1^w}{4}} \Gamma\paren{\tfrac{2-2s_3+\mu_2^w-\mu_1^w}{4}} \Gamma\paren{\tfrac{2-2s_3+\mu_1^w-\mu_2^w}{4}} \frac{ds_3}{2\pi i} \\
	&= \pi \, \Gamma\paren{\tfrac{1+\mu_2^w-\mu_1^w}{2}} \Gamma\paren{\tfrac{1+\mu_1^w-\mu_2^w}{2}},
\end{align*}
so
\begin{align*}
	K_{w_4}(y,\mu) \sim& \frac{1}{64\pi} \sum_{w\in W_3} \abs{\pi^3 y_1}^{1-\mu_3^w} \frac{\Gamma\paren{\tfrac{\mu_3^w-\mu_1^w}{2}} \Gamma\paren{\tfrac{\mu_3^w-\mu_2^w}{2}}}{\Gamma\paren{\tfrac{1+\mu_1^w-\mu_3^w}{2}} \Gamma\paren{\tfrac{1+\mu_2^w-\mu_3^w}{2}}}.
\end{align*}
Thus we conclude
\begin{align}
\label{eq:C3Value}
	C_3(\mu,-1) =& \frac{1}{512} \frac{1}{\sin\frac{\pi}{2}(\mu_1-\mu_3) \sin\frac{\pi}{2}(\mu_2-\mu_3)},
\end{align}
after some simplification.

\section{Mellin-Barnes integrals}
\label{sect:MBRepns}
Theorem \ref{thm:LongEleMBs} may be checked against the power series representations by shifting contours to the left.
Some care must be taken that
\[ \Re(2s_1-s_2)<0, \qquad \Re(2s_2-s_1) < 0, \]
to maintain absolute convergence (see the note on contours following Theorem \ref{thm:LongEleMBs} for the initial contour); this requires shifting the contours in stages.

The integral of Theorem \ref{thm:w4MBs} is a little more complicated:
To obtain this last representation, let $\mathcal{C}$ be the contour which travels along straight lines from $-\infty-iM$ to $\eta-iM$ to $\eta+iM$ to $-\infty+iM$ with $M>\max_i\abs{\Im(\mu_i)}$, then
\begin{align*}
	J_{w_4}(y_1,\mu) =& \frac{1}{\pi^2} \int_\mathcal{C} \abs{8\pi^3 y_1}^{1-s} \exp\paren{i\tfrac{\pi}{2}\sgn(y_1) s} \Gamma(s-\mu_1)\Gamma(s-\mu_2)\Gamma(s-\mu_3) \\
	& \qquad \exp\paren{-i\tfrac{\pi}{2}\sgn(y_1) \mu_3} \sin\pi(s-\mu_1)\sin\pi(s-\mu_2) \frac{ds}{2\pi i}.
\end{align*}
The representation of $K_{w_4}$ follows by trigonometry and straightening the contour.
What allows us to straighten the contour is that the terms involving $\exp\paren{i\tfrac{5\pi}{2}\sgn(y_1) s}$ cancel in the sum over the Weyl group.

\section{Aside on the Eisenstein series constant}
\label{sect:EisenConsts}
We briefly paraphrase \cite[sects. 10.10-10.13]{Goldfeld} to fix the value of the constant on the minimal parabolic Eisenstein series term in Theorem \ref{thm:LisKuznetsovFormula}.
Suppose $f:\Gamma\backslash G/K\to\C$ is orthogonal to the residues of the Eisenstein series, and consider the $1,1,1$ constant term
\[ f^{1,1,1}(y) := \int_{U(\Z)\backslash U(\R)} f(uy) du, \]
and its Mellin expansion
\[ f^{1,1,1}(y) = \frac{1}{(2\pi i)^2} \int_{(2)} \int_{(-2)} y_1^{1-\mu_3} y_2^{1+\mu_1} \int_0^\infty \int_0^\infty f^{1,1,1}(t) t_1^{1+\mu_3} t_2^{1-\mu_1} \frac{dt_1 \, dt_2}{(t_1 t_2)^3} d\mu_1 \, d\mu_3. \]
By folding, and the assumptions on $f$, this is
\[ f^{1,1,1}(y) = \frac{1}{(2\pi i)^2} \int_{\Re(\mu)=0} p_{\rho+\mu}(y) \int_{\Gamma\backslash G/K} f(z) E(z;-\mu) dz \, d\mu. \]
Now the $1,1,1$ constant term of the Eisenstein series is
\[ E^{1,1,1}(z;\mu) = 4 p_{\rho+\mu}(y) + \text{other Weyl element terms}, \]
the 4 here is the number of the diagonal orthogonal matrices, $\abs{V}$.
Its functional equations are the transforms which permute the Weyl element terms of the $1,1,1$ constant term and these transforms are orthogonal on $\Re(\mu)=0$, so we consider
\begin{align*}
	f^{1,1,1}(y) =& \frac{1}{(2\pi i)^2} \frac{1}{24}\paren{\sum_{w\in W} 1} \int_{\Re(\mu)=0} 4 p_{\rho+\mu}(y) \int_{\Gamma\backslash G/K} f(z) E(z;-\mu) dz \, d\mu \\
	=& \frac{1}{24 (2\pi i)^2} \int_{\Re(\mu)=0} E^{1,1,1}(z;\mu) \int_{\Gamma\backslash G/K} f(z) E(z;-\mu) dz \, d\mu,
\end{align*}
after substituting $\mu \mapsto \mu^w$ and applying the functional equations of $E(z;-\mu^w)$.

\section{Acknowledgements}
The author would like to thank Valentin Blomer for his encouragement and many helpful discussions; in particular, for the connection with the Voronoi formula in \eqref{eq:Kw4Voronoi}.

\bibliographystyle{amsplain}

\bibliography{SpectralKuznetsov}

\end{document}